\documentclass[11pt]{amsart}

\usepackage{amssymb}
\usepackage[mathscr]{eucal}

\numberwithin{equation}{section}

\newtheorem{theorem}{Theorem}[section]
\newtheorem*{theorem*}{Theorem}
\newtheorem{proposition}[theorem]{Proposition}
\newtheorem{corollary}[theorem]{Corollary}
\newtheorem{lemma}[theorem]{Lemma}

\newtheorem{definition}{Definition}[section]
\newtheorem{condition}{Condition}[section]

\newcommand{\R}{\mathbb{R}}
\newcommand{\Q}{\mathbb{Q}}
\newcommand{\C}{\mathbb{C}}
\newcommand{\Z}{\mathbb{Z}}
\newcommand{\p}{\partial}
\newcommand{\half}{\tfrac{1}{2}}
\newcommand{\CP}{\mathbb{CP}}

\renewcommand{\L}{\mathcal{L}}
\newcommand{\CL}{\mathcal{L}}
\newcommand{\<}{\langle}
\renewcommand{\>}{\rangle}

\DeclareMathOperator{\Tr}{Tr}
\DeclareMathOperator{\ad}{ad}

\DeclareMathOperator{\vdim}{vdim}

\renewcommand{\o}{\otimes}
\newcommand{\Mbar}{\overline{\mathcal{M}}}
\newcommand{\Mtilde}{\widetilde{\mathcal{M}}}

\newcommand{\PHI}{\mathcal{F}}

\newcommand{\tTheta}{\tilde{\Theta}}
\renewcommand{\t}{t}
\newcommand{\om}{\omega}
\newcommand{\Om}{\Omega}

\newcommand{\virt}{{\textup{virt}}}
\newcommand{\bull}{\bullet}
\newcommand{\CO}{\mathcal{O}}
\newcommand{\Ohat}{\widehat{\mathcal{O}}}
\DeclareMathOperator{\ev}{ev}
\newcommand{\tint}{{\textstyle\int}}

\newcommand{\A}[1]{\mathcal{A}_{{#1}}}
\renewcommand{\AA}{\mathcal{A}}
\newcommand{\Id}{I}
\renewcommand{\]}{{]\!]}}

\DeclareMathOperator{\End}{End}
\DeclareMathOperator{\Spec}{Spec}

\DeclareMathOperator{\Res}{Res}

\newcommand{\stab}{{\textup{stab}}}
\newcommand{\pt}{\mathsf{pt}}
\newcommand{\WW}{\mathcal{W}}
\renewcommand{\.}{\wedge}

\newcommand{\U}{\mathcal{U}}
\newcommand{\Dil}{\mathcal{Q}}

\renewcommand{\P}{\pi}

\newcommand{\g}{\mathbf{t}}
\newcommand{\G}{\mathbf{s}}
\newcommand{\GG}{\mathcal{G}}

\newcommand{\D}{\sigma}
\renewcommand{\d}{\tau}
\newcommand{\eps}{\varepsilon}

\renewcommand{\*}{\cdot}
\newcommand{\tdelta}{\tilde{\delta}}

\newcommand{\CH}{\mathcal{H}}
\renewcommand{\H}{\mathsf{H}}

\newcommand{\Nov}{\Lambda}

\newcommand{\normal}[1]{\textup{\,:\,}#1\textup{\,:\,}}

\newcommand{\YY}{\mathcal{V}}
\newcommand{\Eta}{\pmb{\eta}}
\newcommand{\V}{\mathcal{V}}
\newcommand{\PP}{\mathcal{P}}
\newcommand{\X}{\mathcal{X}}
\newcommand{\M}{\mathcal{M}}
\newcommand{\h}{\mathcal{H}}
\newcommand{\quart}{\tfrac{1}{4}}

\DeclareMathOperator{\Aut}{Aut}
\let\Vert\undef
\DeclareMathOperator{\Vert}{Vert}

\begin{document}

\title[Higher-genus Gromov-Witten invariants]{The jet-space of a
Frobenius manifold and higher-genus Gromov-Witten invariants}

\author{Ezra Getzler}

\address{Dept. of Mathematics, Northwestern University, Evanston, IL
60208, USA}

\email{getzler@math.northwestern.edu}

\dedicatory{To Yuri I. Manin on his birthday}

\maketitle

\section{Introduction}

The theory of genus zero Gromov-Witten invariants associates to a compact
symplectic manifold $X$ a Frobenius manifold $\H$ (also known as the small
phase space of $X$) whose underlying flat manifold is the cohomology space
$H^*(X,\C)$. Higher genus Gromov-Witten invariants give rise to a sequence
of generating functions $\PHI^X_g$, one for each genus $g>0$; these are
functions on the large phase space
$$
\H_\infty = H^*(X\times\CP^\infty,\C) .
$$
The manifold $\H_\infty$ has a rich geometric structure: it is the
jet-space of curves in the Frobenius manifold $\H$. (This identification is
implicit in Dubrovin \cite{D}.)

Frobenius manifolds also arise naturally in a number of other geometric
situations, such as singularity theory and mirror symmetry (see for example
Dubrovin \cite{D} and Manin~\cite{M}). The Frobenius manifolds of
Gromov-Witten theory carry an additional geometric structure, a fundamental
solution, given by the formula
$$
\Theta_a^b(z) = \delta_a^b + \sum_{k=0}^\infty z^{k+1}
\<\<\tau_k(\CO_a)\CO^b\>\>_0^X .
$$
In this paper, we formulate the differential equations satisfied by the
potentials $\PHI^X_g$, such as topological recursion relations and the
Virasoro constraints, in an intrinsic fashion, that is, in such a way
that the equations do not depend on the choice of fundamental
solution. This effort is rewarded by a closer relationship between the
resulting theory and the geometry of moduli spaces of stable curves.

A consequence of our analysis is the proof of a conjecture of Eguchi and
Xiong \cite{EX}. (Here, we use the summation convention with respect to
indices $a$, $b$, \dots.) Introduce the coordinates
$$
u^a = \p \<\<\CO^a\>\>_0^X
$$
on $\H$, along with their derivatives
$$
\p^ku^a = \p^{k+1}\<\<\CO^a\>\>_0^X ,
$$
which form a coordinate system on $\H_\infty$; here, $\p$ is the vector
field on $\H_\infty$ given by differentiation with respect to the puncture
variable $t^e_0$.
\begin{theorem} \label{Fg}
If $g>0$, the Gromov-Witten potential $\PHI^X_g$ has the form
$$
\PHI^X_g = \sum_{n=0}^{3g-3} \frac{1}{n!}
\sum_{\substack{k_i>0\\k_1+\dots+k_n\le3g-3}} f^{k_1\dots k_n}_{a_1\dots
a_n}(u^a,\p u^a) \, \p^{k_1+1}u^{a_1} \dots \p^{k_n+1}u^{a_n} .
$$
The coefficients $f^{k_1\dots k_n}_{a_1\dots a_n}(u^a,\p u^a)$ are
symmetric in the indices $(k_i,a_i)$, and homogeneous of degree
$(2g-2)-(k_1+\dots+k_n+n)$ in the variables $\p u^a$.
\end{theorem}

This theorem generalizes a well-known formula for the case of pure gravity:
by the Kontsevich-Witten theorem, the functions $f^{k_1\dots k_n}_{a_1\dots
a_n}(u,\p u)$ are coefficients of the Gelfand-Dikii polynomials
(Hamiltonians of the KdV hierarchy), and the theorem may be reinterpreted
as a result in the theory of integrable hierarchies. Of course, it may be
easily proved in this case by direct methods from the definition of the
Gelfand-Dikii polynomials. Dubrovin and Zhang \cite{DZ2} have greatly
generalized this observation, showing that Theorem \ref{Fg} holds, for
quite different reasons, in the theory of integrable systems.

\subsection*{Acknowledgements}
Many of the ideas of this paper come from collaboration with
T.~Egu\-chi and C.-S. Xiong, and from conversations with
B.~Du\-bro\-vin, X.~Liu and Y.~Zhang.

The results of the first four sections of this paper were obtained during
the author's residence at RIMS, Kyoto, during the 1999--2000 year devoted
to ``Geometry Related to String Theory''; I thank Kyoji Saito and Masa-Hiko
Saito for the invitation to participate. These results were presented at
conferences at RIMS in June 2000, organized by M.-H. Saito, and at KIAS,
Seoul, in August 2000, organized by Y.-G.~Oh.

The results of the remaining sections were obtained were obtained during a
membership of the Institute for Advanced Study, Princeton, in Winter, 2002.

This work was partially supported by the NSF under grants DMS-9794320 and
DMS-0072508, and, through the IAS, under grant DMS-9729992.

\section{The jet-space of a Dubrovin manifold}

\subsection{Dubrovin connections} \label{Dubrovin}

A large part of the theory of Frobenius manifolds does not require the
existence of a metric or Euler vector field. In this section, we introduce
the geometric structure which is the essential part of the definition of a
Frobenius manifold, a Dubrovin connection.

As in Dubrovin \cite{pencil}, we work with connections on the cotangent
bundle; of course, there is a bijection between connections on the tangent
and cotangent bundles. Let
$$
\nabla : \Gamma(M,T^*M)\to\Om^1(M,T^*M)
$$
be a connection on the cotangent bundle of a manifold $M$. Denote by
$$
\nabla_X = \iota(X)\nabla:\Gamma(M,T^*M)\to\Gamma(M,T^*M)
$$
the operation of covariant differentiation along a vector field $X$ on
$M$.

The connection $\nabla$ is torsion-free if the \textbf{torsion}
$T\in\Om^2(M,TM)$ vanishes, and flat if the \textbf{curvature}
$R\in\Om^2(M,\End(T^*M))$ vanishes. A \textbf{flat} manifold is a manifold
$M$ with torsion-free and flat connection $\nabla$ on the cotangent bundle
$T^*M$. A flat coordinate chart $M\supset U\to\R^N$ on a manifold is a
coordinate chart $u^a:U\to M$ such that the one-forms $du^a$ are parallel:
$$
\nabla(du^b) = 0 .
$$
A manifold is flat if and only if it has an atlas of flat coordinate charts.

A \textbf{pencil} of torsion-free affine connections is a family of
connections $\nabla^z$ depending on a complex parameter $z$ such that
$$
\nabla^z = \nabla + z\AA ,
$$
where $\nabla=\nabla^{z=0}$ and $\AA\in\Om^1(M,\End(T^*M))$. A
\textbf{Dubrovin connection} is a pencil $\nabla^z$ of torsion-free
affine connections which is flat for all values of $z$, together with
a flat vector field $e$ such that $\A{e}$ is the identity. A Dubrovin
manifold $(M,\nabla^z,e)$ is a manifold together with a Dubrovin
connection.

Given a Dubrovin connection $\nabla^z$, denote the vector field
$(\A{X})^*Y$ by $X\circ Y$; this defines a product on the tangent bundle of
$M$. Denote the components of the tensor $\AA$ in flat coordinates $(u^a)$
by $\AA_{ab}^c$:
$$
\p_a\circ\p_b = \AA_{ab}^c\p_c .
$$
\begin{proposition} \label{dubrovin}
The data $(\nabla^z,e)$ form a Dubrovin connection if and only if the
following conditions are satisfied:
\begin{enumerate}
\item the affine connection $\nabla$ is flat;
\item the product $X\circ Y$ is commutative and associative:
$$
(X\circ Y) \circ Z = X\circ(Y\circ Z) ;
$$
\item the vector field $e$ is flat, $\nabla e=0$, and is an identity for
the product $X\circ Y$:
$$
e\circ X = X\circ e = X ;
$$
\item $[\nabla_X,\A{Y}]-[\nabla_Y,\A{X}]=\A{[X,Y]}$.
\end{enumerate}
\end{proposition}
\begin{proof}
If $\nabla$ is torsion-free, the commutativity of the product $X\circ
Y$ is equivalent to the vanishing of the torsion of $\nabla^z$.

The curvature $R^z(X,Y)$ of the connection $\nabla^z$ is a quadratic
polynomial in $z$:
\begin{multline*}
[\nabla^z_X,\nabla^z_Y]-[\nabla^z_X,\nabla^z_Y]-\nabla^z_{[X,Y]}
\\ = R(X,Y) + z\bigl( [\nabla_X,\A{Y}]-[\nabla_Y,\A{X}]-\A{[X,Y]} \bigr) + z^2
[\A{X},\A{Y}] .
\end{multline*}
The constant term vanishes if and only if $\nabla$ is flat, the linear
term if and only if 4) holds, and the quadratic term if and
only if $X\circ Y$ is associative.
\end{proof}

The questions which we address in this paper are of local character; for
this reason, we will work in the neighbourhood of a basepoint $p\in M$ of
our Dubrovin manifold $(M,\nabla^z,e)$. Let $(u^a)$ be a flat coordinate
system on $M$ centred at $p$: that is, $u^a(p)=0$.

Since the identity vector field $e$ is flat, we may assume that the flat
coordinate system $(u^a)$ is chosen in such a way that $e$ is
differentiation with respect to $u^e$, where $e$ is one of the indices
labelling the coordinates: in other words, $e=\p/\p u_e=\p_e$.

The one-form
\begin{equation} \label{socle}
\om = \Tr(\A{a})du^a\in\Om^1(M)
\end{equation}
plays a prominent role in the theory of Dubrovin manifolds; following
Hert\-ling \cite{Hertling}, we call it the \textbf{socle} one-form.

\subsection{Fundamental solutions}
A fundamental solution of a Dubrovin manifold is a power series
$$
\Theta = \Id + \sum_{n=0}^\infty z^{n+1} \Theta_n
\in\Gamma(M,\End(T^*M))\[z\]
$$
such that
$$
\nabla \* \Theta(z) = \Theta(z) \* \nabla^z .
$$
This is equivalent to the sequence of equations
\begin{equation} \label{Thetan}
\nabla_X \Theta_n = \Theta_{n-1} \A{X} .
\end{equation}
Since the endomorphism $\Theta_0$ plays a special role, we introduce the
notation $\M=\Theta_0$. The case $n=0$ of \eqref{Thetan} shows that
$\nabla_X\M=\A{X}$. In particular,
$$
\p_a \M_e^b = \delta_a^b .
$$
It follows that $\M_e^a=u^a+c^a$, where $c^a$ is a constant.
\begin{proposition} \label{Fundamental}
Let $p$ be a point of a Dubrovin manifold $(M,\nabla^z,e)$, and let $(u^a)$
be a system of flat coordinates which vanish at $p$. Then there exists a
fundamental solution $\Theta(z)$ in a neighbourhood of $p$ such that
$\M_e^a=u^a$.
\end{proposition}
\begin{proof}
We construct endomorphisms $\Theta_n$ of the cotangent bundle $T^*M$
inductively in $n$, starting with $\Theta_{-1}=\Id$. To carry out the
induction, we must show that the one-form
$\Theta_{n-1}\AA\in\Om^1(M,\End(T^*M))$ is exact. To do this, we use
condition 4) of Proposition \ref{dubrovin}: we have
\begin{align*}
\nabla\bigl( \Theta_{n-1}\AA \bigr) &= \bigl( [\nabla_a,\Theta_{n-1}\A{b}]
- [\nabla_b,\Theta_{n-1}\A{a}] \bigr) \, du^a \. du^b \\
&= \bigl( \Theta_{n-1} \A{[\p_a,\p_b]} +
\Theta_{n-2} [\A{a},\A{b}] \bigr) \, du^a \. du^b = 0 .
\end{align*}
In particular, the existence of $\Theta_0$ is equivalent to condition 4) of
Proposition \ref{dubrovin}, while granted the existence of $\Theta_0$, the
existence of $\Theta_1$ is equivalent to the associativity of the product
$X\circ Y$.
\end{proof}

If $\Theta(z)$ and $\tTheta(z)$ are two fundamental solutions of the
Dubrovin connection $\nabla^z$, then $\tTheta(z)^{-1}\Theta(z)$ is a flat
section of $\Gamma(M,\End(T^*M))\[z\]$, that is,
$$
\nabla(\tTheta(z)^{-1}\Theta(z)) = 0 .
$$
Conversely, if $\rho(z)$ is a flat section of
$\Gamma(M,\End(T^*M))\[z\]$ such that $\rho(0)=\Id$ and $\Theta(z)$ is a
fundamental solution, then $\Theta(z)\rho(z)$ is again a fundamental
solution.

\subsection{The jet-space of a Dubrovin manifold}

Let $n$ be a natural number and let $M$ be a manifold. An $n$-jet in $M$ is
a map from the variety $\Spec(\C[t]/(t^{n+1}))$ to $M$.

In local coordinates $(u^a)$, we may write an $n$-jet as
$$
u^a(t) = \sum_{k=0}^n \frac{t^k}{k!} u^a_k + O(t^{n+1}) .
$$
In particular, $u^a_0$ are the coordinates of the origin of the jet, and
$u^a_1\,\p_a\in T_{u(0)}M$ is the velocity of the jet at $t=0$.

The space of $n$-jets $J^nM$ in $M$ is a fibre bundle over $M$ whose
fibres are affine spaces of dimension $n\dim M$.  For example, a
$0$-jet is the same as a point of $M$, and a $1$-jet is the same as a
tangent vector; thus $J^0M=M$ and $J^1M=TM$. If $m<n$, denote the
projection from $J^nM$ to $J^mM$ defined by reduction mod $t^{m+1}$ by
$\rho_{n,m}$; we write $\rho_n:J^nM\to M$ instead of $\rho_{n,0}$.

Let $J^\infty M$ be the inverse limit
$$
\varprojlim_n J^nM .
$$
There is a fibration $\rho_{\infty,n}$ from $J^\infty M$ to $J^nM$.

If $\V$ is a vector bundle on $M$, denote the space of sections of the
vector bundle $\rho_n^*\V$ on $J^nM$ by $\V_n$; in particular, the space of
sections
$$
\V_\infty = \Gamma(J^\infty M,\rho_\infty^*\V) = \bigcup_{n=0}^\infty \V_n
$$
is filtered by subspaces
$$
\V_0 = \Gamma(M,\V) \subset \V_1 \subset \dots \subset \V_n \dots \subset
\V_\infty .
$$
In the case where $\V$ is the trivial bundle $\CO$, we see that $\CO_n$ is
the algebra of functions on $J^nM$.

Given a coordinate system $(u^a)$ on $M$, define vector fields
$$
\p_{\ell,a} = \frac{\p}{\p u^a_\ell} , \quad 0\le\ell\le n ,
$$
on the jet-space $J^nM$. The differential operator
$$
\p = \sum_{\ell=0}^\infty u^a_{\ell+1} \, \p_{\ell,a}
$$
is a vector field on $J^\infty M$; formally speaking, $\p$ generates
translation along the parameter $t$ of the jet. Note that $u^a_k=\p^ku^a$,
where $u^a=u^a_0$.

Another important vector field on the jet-space $J^\infty M$ is the
\textbf{dilaton} vector field
$$
\Dil = \sum_{\ell=0}^\infty \ell \, u^a_\ell \, \p_{\ell,a} ,
$$
which counts the coordinate $u^a_\ell=\p^\ell u^a$ with weight
$\ell$.

We now specialize to the case where $M$ is a Dubrovin manifold. The
Dubrovin connection $\nabla^z$ on $T^*M$ pulls back by $\rho_\infty$ to a
connection $\nabla^z$ on $\rho_\infty^*T^*M$. If $X$ is a vector field on
the jet-space $J^\infty M$, denote by $X$ the covariant derivative
$\nabla_X$ on $\rho_\infty^*T^*M$ and its associated tensor bundles. For
example, we write $\p$ instead of $\nabla_\p$.

Given a vector field $X$ on $J^\infty M$, define an endomorphism $\A{X}$
of the bundle $\rho_\infty^*\End(T^*M)$ by the formula
$$
\A{X} = X(\M) .
$$
This definition is independent of the fundamental solution $\Theta(z)$,
and consistent with the definition of $\A{X}$ in the case that $X$ is a
vector field on $M$. We have the formula
$$
X(\Theta(z)) = z\Theta(z)\A{X} \in \End(T^*M)_\infty .
$$

If $M$ has basepoint $p$ and $n>0$, let $\P$ be the basepoint of $J^nM$
such that $\P(t)=p+te+O(t^{n+1})$. At $\P$, the endomorphism
$\X=\A{\p}\in\End(T^*M)_\infty$ equals the identity; it follows that for
$n>0$, the open subset $J^n_0M\subset J^nM$ on which $\X$ is invertible is
a neighbourhood of $\P\in J^nM$. If $\V$ is a vector bundle on $M$, denote
the space of sections of $\rho_n^*\V$ over $J^n_0M$ by $\V_{n,0}$. In
particular, we have
$$
\V_{\infty,0} = \bigcup_{n=0}^\infty \V_{n,0} .
$$

\subsection{A frame of vector fields on the jet-space}
Let $(M,\nabla^z,e)$ be a Dubrovin manifold with basepoint $p$, and let
$\Theta(z)$ be a fundamental solution such that the flat coordinates
$\M_e^a=u^a$ vanish at $p$. Consider the following noncommutative analogues
of the Fa\`a di Bruno polynomials:
$$
\PP_n(z) = z^{-1} \Theta^{-1}(z)\p^n\Theta(z) \in \End(T^*M)_n .
$$
For example, we have $\PP_0(z)=z^{-1}\Id$ and $\PP_1(z)=\X$. In general, we
see that for $n>0$,
$$
\PP_n(z) = (\p+z\X)^{n-1}\X = \sum_{j=0}^{n-1} z^j \PP_{n,j} ,
$$
where $\PP_{n,j}\in\End(T^*M)_{n-j}$. The polynomials $\PP_n(z)$ are
intrinsic; they do not depend on the choice of fundamental solution.

Let $\p_\ell = du^a \, \p_{\ell,a}$. We now introduce a generating
function for vector fields on $J^\infty M$, by the formula
$$
\D(z) = \sum_{k=0}^\infty z^k \D_k = \sum_{\ell=0}^\infty
\PP_{\ell+1}(z) \p_\ell : \CO_\infty \to T^*M\[z\]_\infty .
$$
For example,
\begin{align*}
\D_0 &= \sum_{\ell=0}^\infty \p^\ell\X \, \p_\ell , &
\D_1 &= \sum_{\ell=1}^\infty \sum_{k=1}^\ell
\p^{\ell-k}(\X\,\p^{k-1}\X) \, \p_\ell .
\end{align*}
Observe that $[\p+z\X,\D(z)]=0$; it follows that
$[\p,\D_k]=-\X\D_{k-1}$. Together with the formula
$$
\D_a(z) u^b = \X_a^b ,
$$
this characterizes $\D(z)$.

If $X$ is vector field on $M$, define vector field $\D_{k,X}$ on
$J^\infty M$ by the generating function
$$
\D_X(z) = \< X , \D(z) \> = \sum_{k=0}^\infty z^k \, \D_{k,X} .
$$
We abbreviate the vector fields $\D_{k,\p_a}$ to $\D_{k,a}$. In
particular, $\D_{0,e}=\p$.

Given of sequences $k=(k_1,\dots,k_n)$ and $\ell=(\ell_1,\dots,\ell_n)$, we
say that $k\ge\ell$ if $k_i\ge\ell_i$ for all $i$, and that $k>\ell$ if in
addition $k_i>\ell_i$ for at least one $i$.
\begin{proposition} \label{vanish}
A function $f\in\CO_\infty=\CO_{\infty,0}$ lies in $\CO_{\ell,0}$ if and
only if $\D(z)f$ is a polynomial in $z$ of degree at most $\ell$.

More generally, suppose $K=(K_1,\dots,K_n)$ is a sequence of integers such
that $K_i>1$ for all $i$, and suppose that $\D_{k_1}\dots\D_{k_n}f=0$ for
all $k\ge K$; then $\p_{k_1}\dots\p_{k_n}f=0$ for all $k\ge K$.
\end{proposition}
\begin{proof}
Since $\D_j\X=0$ for $j>1$, we see that for $k\ge K$,
\begin{multline*}
\D_{k_1,a_1}\dots\D_{k_n,a_n} \\ = \bigl(\X^{k_1+1}\bigr){}_{a_1}^{b_1}
\dots \bigl(\X^{k_n+1}\bigr){}_{a_n}^{b_n}
\p_{k_1,b_1}\dots\p_{k_n,b_n} + \sum_{\ell>k} a(k,\ell) \,
\p_{\ell_1}\dots\p_{\ell_n} .
\end{multline*}
Since $\X$ is invertible on $J^\infty_0M$, it follows that
\begin{multline*}
\p_{k_1,a_1}\dots\p_{k_n,a_n} \\ = \bigl(\X^{-k_1-1}\bigr){}_{a_1}^{b_1}
\dots \bigl(\X^{-k_n-1}\bigr){}_{a_n}^{b_n}
\D_{k_1,b_1}\dots\D_{k_n,b_n} + \sum_{\ell>k} a(k,\ell) \,
\D_{\ell_1}\dots\D_{\ell_n} .
\end{multline*}
Thus, $\D_{k_1}\dots\D_{k_n}f=0$ for $k\ge K$ if and only if
$\p_{k_1}\dots\p_{k_n}f=0$ for $k\ge K$.
\end{proof}

\subsection{An affine structure on the jet-space}
We now show how the choice of a fundamental solution $\Theta(z)$ on a
Dubrovin manifold $(M,\nabla^z,e)$ gives rise to an affine structure on the
jet-space $J^\infty_0M$, that is, a frame of the tangent bundle
$TJ^\infty_0M$ consisting of commuting vector fields. This affine structure
identifies the formal neighbourhood of $\p$ with the large phase space of
Gromov-Witten theory.

Introduce the generating function of vector fields on $J^\infty M$,
$$
\d(z) = \sum_{k=0}^\infty z^k \d_k = \Theta(z)\D(z) = z^{-1}
\sum_{\ell=0}^\infty \p^{\ell+1}\Theta(z) \, \p_\ell .
$$
The introduction of these vector fields is motivated by the fact that they
are evolutionary:
$$
[\p,\d(z)] = 0 .
$$
\begin{lemma} \label{dild}
$[ \Dil , \D(z) ] = \D(z)$ and $[ \Dil , \d(z) ] = \d(z)$
\end{lemma}
\begin{proof}
Since $\Dil\Theta(z)=0$, it suffices to prove that
$[\Dil,\D(z)]=\D(z)$. Since $[\Dil,\p]=\p$, it follows that
$\Dil\PP_\ell=\ell\PP_\ell$, and hence that
\begin{equation*}
[ \Dil , \D(z) ] = \sum_{\ell=0}^\infty [ \Dil , \PP_{\ell+1}(z) \p_\ell ]
= \sum_{\ell=0}^\infty \bigl( (\ell+1)\PP_{\ell+1}(z) \p_\ell - \ell
\PP_{\ell+1}(z) \p_\ell \bigr) = \D(z) .
\qedhere\end{equation*}
\end{proof}

If $X$ is a vector field on $M$, let $\d_{k,X}$ be the vector field on
$J^\infty M$ defined by the generating function
$$
\d_X(z) = \< X , \d(z) \> = \sum_{k=0}^\infty z^k \, \d_{k,X} ,
$$
and let $\d_{k,a}=\d_{k,\p_a}$. Note that $\d_0=\D_0$; in particular,
$\d_{0,e}=\p$. Also, observe that if $f\in\CO_0$ is a function only of
the coordinates $u^a$ on the Frobenius manifold $M$, then
$$
\d_{0,a}f = \D_{0,a}f = \X_a^b \p_b f .
$$
In particular, we see that
\begin{equation} \label{dp}
\p_af = (\X^{-1})_a^b \d_{0,b}f .
\end{equation}
\begin{proposition}
The vector fields $\d_{k,a}$ form a frame of the tangent bundle of
$J^\infty_0M$.
\end{proposition}
\begin{proof}
We have
$$
\D_{k,a} = \bigl(\X^{k+1}\bigr){}_a^b \p_{k,b} + \sum_{\ell>k} a(k,\ell) \,
\p_\ell ,
$$
and the result follows, since $\X$ is invertible on $J^\infty_0M$.
\end{proof}

The vector fields $\d_{k,a}$ were shown by Dubrovin \cite{D} to be the
commuting flows of an integrable hierarchy. In this paper, these vector
fields play an auxilliary role: the vector fields $\D_{k,a}$ are a more
natural frame for the tangent bundle of $J^\infty_0M$, because they are
intrinsic, that is, independent of the fundamental solution
$\Theta(z)$. For us, the fact the vector fields $\D_{k,a}$ do not commute
is of secondary importance.

\begin{proposition} \label{involution}
$[\d_a(z),\d_b(y)]=0$
\end{proposition}
\begin{proof}
We have $\A{\D_a(z)}=\A{\D_{0,a}}=\WW_a=\X\A{a}$; in other words, the
components $\WW_{ab}^c$ of $\WW_a$ equal $\AA_{ab}^{c'}\X_{c'}^c$, and are
symmetric in $a$ and $b$. Since $[\p,\d(z)]=0$, we see that
\begin{align*}
[\d_a(y),\d_b(z)] &= \sum_{k=0}^\infty \p^{k+1} \bigl( z^{-1}
\d_a(y)\Theta_b^c(z) - y^{-1} \d_b(z)\Theta_a^c(y) \bigr) \p_{k,c}
\\
&= \sum_{k=0}^\infty \p^{k+1} \bigl( z^{-1}
\Theta_a^{a'}(y)\D_{a'}(y)\Theta_b^c(z) - y^{-1}
\Theta_b^{b'}(z)\D_b(z)\Theta_a^c(y) \bigr) \p_{k,c} \\
&= \sum_{k=0}^\infty \p^{k+1} \bigl( \Theta_a^{a'}(y)\Theta_b^{b'}(z)
\bigl( \WW_{a'b'}^c - \WW_{b'a'}^c \bigr) \bigr) \p_{k,c} = 0 .
\qedhere\end{align*}
\end{proof}

\begin{corollary}
We have $[\D_a(y),\D_b(z)]=\WW_{ab}^c(y\D_c(y)-z\D_c(z))$; in other words,
$$
[\D_{k,a},\D_{\ell,b}] = \begin{cases}
- \WW_{ab}^c \D_{\ell-1,c} , & k=0 , \\
0 , & k>0 . \end{cases}
$$
\end{corollary}
\begin{proof}
Since
$[\d_a(y),\d_b(z)]=[\Theta_a^{a'}(y)\D_{a'}(y),\Theta_b^{b'}(z)\D_{b'}(z)]$
vanishes, we see that
\begin{multline*}
\Theta_a^{a'}(y)\Theta_b^{b'}(z)[\D_{a'}(y),\D_{b'}(z)] \\
\begin{aligned}
&= \Theta_b^{b'}(z)\bigl(\D_{b'}(z)\Theta_a^{a'}(y)\bigr)\D_{a'}(y) -
\Theta_a^{a'}(y) \bigl( \D_{a'}(y)\Theta_b^{b'}(z) \bigr) \D_{b'}(z) \\
&= \Theta_a^{a'}(y)\Theta_b^{b'}(z) \bigl( y\WW_{b'a'}^c\D_c(y) -
z\WW_{a'b'}^c\D_c(z) \bigr) .
\qedhere\end{aligned}
\end{multline*}
\end{proof}

\subsection{Flat coordinates on the jet-space}
The results of this section are taken from \cite{virasoro}. The main result
of this section is Theorem \ref{G}, which is taken from Section 4 of
\cite{virasoro}; it is interesting to compare this with Proposition 6.1 of
Coates and Givental \cite{CG}.

Since the vector fields $\d_{k,a}$ are in involution, there is a coordinate
system $\t_k^a$ centered at $\P\in J^\infty_0M$ such that
$\d_{k,a}=\p/\p\t_k^a$ and $\t_k^a(\P)=0$; we call these the flat
coordinates on the jet-space. Of course, these coordinates do not lie in
$\CO_{\infty,0}$; rather, they lie in the completion $\Ohat$ of
$\CO_\infty$ at $\P$, and are only coordinates in a formal neighbourhood of
$\P$. This completion $\Ohat$ may be identified with the algebra
$\C\[t_k^a\mid k\ge0\]$.

\begin{proposition}
Let $\Theta(z)$ be a fundamental solution on the Dubrovin manifold
$(M,\nabla^z,e)$. There is a unique generating function
$$
\g(z) = \sum_{k=-\infty}^\infty z^{-k-1} \, \g_k \in
TM\[z,z^{-1}\]_{\infty,0}
$$
such that $\Dil\g(z)=-\g(z)$ and
\begin{equation} \label{dg}
\d(y)\g(z) = \frac{\Theta(y)\Theta^{-1}(z)}{z-y}
= \sum_{k=0}^\infty y^k z^{-k-1} \, \Theta(y)\Theta^{-1}(z) .
\end{equation}
For $k\ge0$, we have $\g_k^a=t_k^a-\delta_{k,1}\delta_e^a$.
\end{proposition}
\begin{proof}
We start with a lemma of Dubrovin \cite{D}.
\begin{lemma} \label{integrable}
Let
$$
\Om(y,z) = \sum_{k,\ell=0}^\infty y^k z^\ell \Om_{k,\ell}
$$
be the generating function characterized by the equation
\begin{equation} \label{Om}
(y-z) \Om(y,z) = \Theta(y)\Theta^{-1}(z) - \Id .
\end{equation}
Then $\d_{k,a}\Om_{\ell,m|b}^c = \d_{\ell,b}\Om_{k,m|a}^c$.
\end{lemma}
\begin{proof}
We must prove that $\d_a(x)\Om_b^c(y,z) = \d_b(y)\Om_a^c(x,z)$ or
equivalently, that
$$
(x-z) \d_a(x)\bigl(\Theta(y)\Theta^{-1}(z)\bigr){}_b^c
= (y-z) \d_b(y)\bigl(\Theta(x)\Theta^{-1}(z)\bigr){}_a^c .
$$
But
$$
\d_a(x)\bigl(\Theta(y)\Theta^{-1}(z)\bigr) = (y-z)
\Theta(y) \A{\d_a(x)} \Theta^{-1}(z) .
$$
Thus, it suffices to show that $\bigl(\Theta(y)\A{\d_a(x)}\bigr){}_b^c =
\bigl(\Theta(x)\A{\d_b(y)}\bigr){}_a^c$; this follows from the formula
$$
\bigl( \Theta(y) \A{\d_a(x)} \bigr){}_b^c = \Theta_a^{a'}(x)
\Theta_b^{b'}(y) \A{\D_{0,a'}}{}_{b'}^c = \Theta_a^{a'}(x)
\Theta_b^{b'}(y) \WW_{a'b'}^c
$$
and the symmetry $\WW_{ab}^c=\WW_{ba}^c$.
\end{proof}

To show the existence of $\g(z)$, we must solve the system of equations
\begin{equation} \label{dt}
\d_{j,a} \g_k^b = \delta_{j,k} \delta_a^b - \Om_{j,-k-1|a}^b ,
\end{equation}
subject to the constraints $\Dil\g_k=-\g_k$. For $k\ge0$, this
system has the unique solution
$$
\g_k^a = \t_k^a - \delta_{k,1}\delta_e^a .
$$
For $k<0$, the system is integrable by Lemma~\ref{integrable}. Since
$[\Dil,\d(z)]=\d(z)$ and $\Dil\Om(y,z)=0$, we see that $\Dil\g_k+\g_k$ is
a constant; replacing $\g_k$ by $-\Dil\g_k$, we obtain the desired solution.
\end{proof}

The vector fields $\Dil$ and $e$ have simple expressions in flat
coordinates on the large phase space.
\begin{proposition} \label{flat}
The dilaton vector field $\Dil$ and the identity vector field $e$ are the
residues of the generating functions $-\<\g(z),\d(z)\>$ and
$-\<\g(z),z\d(z)\>$ respectively.
\end{proposition}
\begin{proof}
By Lemma \ref{dild}, there are constants $c_k^a$ such that
$$
\Dil = \sum_{k=0}^\infty \bigl( c_k^a - \t_k^a \bigr) \, \d_{k,a} .
$$
By the definition of $\d(z)$, we have
$$
\d_{k,a} u^b_n = \bigl( \Theta_{k-n-1}\X^{n+1} \bigr){}_a^b .
$$
Applying the equation for $\Dil$ to the function $u^b_n$ and evaluating
at $\P$, we see that
$$
\delta_{n,1} \delta_e^b = \sum_{k=0}^\infty
\Theta_{k-n-1,a}^b(\P) \, c_k^a .
$$
This system of linear equations is upper triangular, with diagonal entries
equal to $1$, hence has a unique solution: it is easily checked that this
solution is $c_k^a=\delta_{k,1}\delta_e^a$, since
$\Theta_e^a(\P)=u^a(\P)=0$.

The proof of the formula for $e$ in flat coordinates is similar. Since
$e(\Theta(z))=z\Theta(z)$ and $[\p,e]=0$, we see that
$$
[e,\d(z)] = z \sum_{\ell=0}^\infty \p^{\ell+1} e(\Theta(z))
\p_\ell = z \d(z) .
$$
In other words, $[e,\d_{k,a}]=\d_{k-1,a}$. It follows that there are
constants $c_k^a$ such that
$$
e = \sum_{k=1}^\infty \bigl( c_k^a - \t_k^a \bigr) \, \d_{k-1,a} .
$$
Applying this equation to the function $u^b_n$ and evaluating at
$\P$, we see that
$$
\delta_{n,0} \delta_e^b = \sum_{k=1}^\infty \Theta_{k-n-2,a}^b(\P) \,
c_k^a .
$$
Again, this system has a unique solution $c_k^a=\delta_{k,1}\delta_e^b$.
\end{proof}

Just as the modification $\D=\Theta(z)^{-1}\d(z)$ of $\d(z)$ is intrinsic,
so the modification $\G=\Theta^*(z)\g(z)$ of $\g(z)$ is intrinsic. This
modification was introduced in \cite{virasoro}, where it was denoted
$z^{-1}G^*(-z^{-1})$.
\begin{theorem} \label{G}
Let $(M,\nabla^z,e)$ be a Dubrovin manifold. The generating function
$$
\G(z) = \Theta^*(z)\g(z) = \sum_{k=-\infty}^\infty z^{-k-1} \G_k \in
TM\[z^{-1}\]_{\infty,0}
$$
is intrinsic, that is, independent of the fundamental solution $\Theta(z)$,
and
$$
\G_k = \begin{cases}
-\<e,(\X^{-1}\p)^{k-1}\X^{-1}\> , & k>0 , \\
0 , & k\le0 .
\end{cases}$$
\end{theorem}
\begin{proof}
We first show that $\G_0=0$. We have
$$
\G_0^a = t_0^a + \sum_{k=1}^\infty t_k^b \Theta_{k-1,b}^a - \Theta_{0,e}^a .
$$
Since $\Om(y,0)=y^{-1}(\Theta(y)-I)$, we see that
$\Om_{k,0}=\Theta_k$, and hence that $\d_k\g_{-1}=-\Theta_k$. It follows
that
$$
\G_0^a = t_0^a + \biggl( \p - \sum_{k=1}^\infty t_k^b \d_{k-1,b} \biggr)
\g_{-1}^a = \g_0^a + e\bigl(\g_{-1}^a\bigr) = 0 .
$$

A special case of \eqref{dg} is $\p\g(z) = z^{-1} \< e , \Theta^{-1}(z)\>$;
it follows that
\begin{equation} \label{GG}
\p\G(z) = z^{-1} e + \X^*\G(z) .
\end{equation}
Taking the residue, we see that $0=e+\X^*\G_1$, while taking the
coefficient of $z^{-k-1}$, we see that $\p\G_k=\X^*\G_{k+1}$ for $k>1$;
this establishes the formula for $\G_k$, $k\ge0$.

In order to prove that $\G_k$ vanishes for $k<0$, we need two lemmas.
\begin{lemma} \label{basic}
A power series $f$ in the coordinates $t_k^a$ such that $\p f=ef=0$ is
constant.
\end{lemma}
\begin{proof}
Since
$$
(\p-e)f = \sum_{k=1}^\infty t_k^a \d_{k-1,a}f
$$
vanishes, it follows that $f$ is a constant. (See \cite{virasoro}, Section
3, for more details of the proof.)
\end{proof}

\begin{lemma}
$(e+z)\g(z)=0$
\end{lemma}
\begin{proof}
The equations $[e,\d(y)]=y\d(y)$ and $(e+z-y)(\Theta(y)\Theta^{-1}(z))=0$
show that
\begin{align*}
(z-y) \d(y)(e+z)\g(z) &= (z-y) (e+z-y)(\d(y)\g(z)) \\
&= (e+z-y)(\Theta(y)\Theta^{-1}(z)) = 0 .
\end{align*}
In other words, $\d(y)(e+z)\g(z)=0$, and hence
$(e+z)\g(z)\in\C\[z,z^{-1}\]$. Composing this equation with the dilaton
vector field $\Dil$, we see that $(e+z)\g(z)=0$.
\end{proof}

It follows from the equations $\Dil\g(z)=-\g(z)$ and $(e+z)\g(z)=0$ that
have $\Dil\G(z)=-\G(z)$ and $e\G(z)=0$.

Suppose that $\G_{-k}$ vanishes, for $k>0$. Taking the coefficient of
$z^{k+1}$ in \eqref{GG}, we see that $\p\G_{-k-1}=0$. Since we know that
$e\G_{-k-1}=0$, it follows from Lemma \ref{basic} that $\G_{-k-1}$ is
constant. Since $\Dil\G_{-k-1}=-\G_{-k-1}$, we conclude that $\G_{-k-1}$
vanishes. Thus $\G_{-k}=0$ for all $k>0$, by induction on $k$.
\end{proof}

Using the series $\G(z)$, we may now rewrite Proposition \ref{flat} in an
intrinsic fashion:
\begin{align*}
\Dil &= - \sum_{k=0}^\infty \< \G_k , \D_k \> , &
e &= - \sum_{k=0}^\infty \< \G_{k+1} , \D_k \> .
\end{align*}
Apart from being intrinsic, these formulas have the virtue that when
applied to a function $f\in\CO_n$, they truncate to a sum over $0\le k\le n$.

\subsection*{Relationship to prior results}
The Dubrovin connection was introduced by Dubrovin \cite{D}. He also
initiated the study of the jet-space of a Frobenius manifold (which he
calls the loop-space). Sections~2.4 and 2.5 are taken from our joint
work with Eguchi and Xiong \cite{EGX}, and Section 2.6 from
\cite{virasoro}, Section 4. (Note that the generating function $\g(z)$
used here is related to the generating function $\theta(\zeta)$ of
\cite{virasoro} by the formula $\g(z)=z^{-1}\theta^*(-z^{-1})$.)

\section{The jet-space of a conformal Dubrovin manifold}

\subsection{Conformal Dubrovin manifolds}
A linear vector field on a flat manifold is a vector field $E$ such that
$\nabla E\in\Gamma(M,\End(TM))$ is flat. In flat coordinates $(u^a)$, there
are constants $A_b^a$ and $B^a$ such that
$$
E = (A_b^au^b+B^a) \p_a .
$$
An \textbf{Euler} vector field $E$ for a Dubrovin connection is
a linear vector field such that
\begin{equation} \label{E}
[E,X\circ Y] = [E,X]\circ Y + X\circ[E,Y] + X\circ Y .
\end{equation}
Equivalently, $[\CL_E,\A{X}]=\A{[E,X]}+\A{X}$, where $\CL_E$ is the Lie
derivative of the vector field $E$ acting on one-forms.

A \textbf{conformal} Dubrovin manifold $(M,\nabla^z,e,E,r)$ is a Dubrovin
manifold $(M,\nabla^z,e)$ together with an Euler vector field $E$ and a
real number $r$.

Let $\mu$ be the endomorphism of the cotangent bundle defined by the
formula
$$
\mu = 1 - \tfrac{r}{2} + \nabla_E - \CL_E .
$$
The adjoint of $\mu$ is the endomorphism of the tangent bundle given
by the formula
$$
\mu^*(X) = - \nabla_EX + [E,X] + ( 1 - \tfrac{r}{2} ) X = ( 1 -
\tfrac{r}{2} )X - \nabla_XE .
$$
Let $\U$ be the endomorphism $\U=\A{E}$ of the cotangent bundle.
\begin{proposition} \label{commute}
Let $(M,\nabla^z,e,E,r)$ be a conformal Dubrovin manifold, and let $\CH$ be
the bundle $T^*M[z,z^{-1}]$. Let $\delta_z$ be the endomorphism
$$
\delta_z = \p_z + z^{-1} \bigl( \mu + \half \bigr) + \U
$$
of $\CH$. Then $[\nabla^z,\delta_z]=0$.
\end{proposition}
\begin{proof}
If $X$ is a vector field on $M$, we have
\begin{multline*}
[\nabla_X^z,\delta_z] = [ \nabla_X + z\A{X} , \p_z + z^{-1} \mu + \A{E} ]
\\ = z^{-1} [\nabla_X,\mu] + \bigl( [\nabla_X,\A{E}] - [\nabla_E,\A{X}] -
\A{X} + [\CL_E , \A{X} ] \bigr) + z [\A{X},\A{E}] .
\end{multline*}
Since $E$ is linear, it follows that $[\nabla_X,\mu]=0$; it is also clear
that $[\A{E},\A{X}]=0$. By \eqref{E}, we see that that
\begin{equation*}
[\CL_E,\A{X}] = \A{[E,X]} + \A{X} = [ \nabla_E , \A{X} ] -
[\nabla_X,\A{E}] + \A{X} .
\qedhere\end{equation*}
\end{proof}

Let $\Theta(z)$ be a fundamental solution of the conformal Dubrovin
manifold $(M,\nabla^z,e,E,r)$, and consider the conjugate of $\delta_z$ by
$\Theta(z)$, defined by the formula
$$
\tdelta_z = \Theta(z) \* \delta_z \* \Theta^{-1}(z) = \p_z + z^{-1} R(z) ,
$$
where
$$
R(z) = \sum_{k=0}^\infty z^k R_k \in z \Gamma(M,\End(T^*M))\[z\] .
$$
In particular, $R_0=\mu+\half$. It follows from Proposition \ref{commute}
that $R(z)$ is flat.

By the Jordan decomposition, we may write $\mu$ as the sum of semisimple
and nilpotent endomorphisms $\mu_s$ and $\mu_n$ of the cotangent bundle
$T^*M$; furthermore, both $\mu_s$ and $\mu_n$ are flat. (In many cases,
such as the theory of Gromov-Witten invariants, $\mu$ is semisimple.) Let
$\tau_\lambda$ be the projection onto the subbundle of $T^*M$ on which
$\mu_s$ has eigenvalue $\lambda$: we have
\begin{align*}
\Id &= \sum_\lambda \pi_\lambda , &
\mu_s &= \sum_\lambda \lambda \pi_\lambda .
\end{align*}
\begin{proposition}
There exists a fundamental solution such that $[\mu_s,R_k]=kR_k$.
\end{proposition}
\begin{proof}
If $\lambda-\lambda'\ne k$, the endomorphism $k-\ad(\mu)$ is invertible on
the space of endomorphisms of the form $\pi_\lambda A \pi_{\lambda'}$: we
have
\begin{align*}
(k-\ad(\mu))^{-1} ( \pi_\lambda A \pi_{\lambda'} ) &= \sum_{i=0}^\infty
(-\ad(\mu_n))^i (k-\ad(\mu_s) ( \pi_\lambda A \pi_{\lambda'} ) )^{-i-1} \\
&= \sum_{i=0}^\infty \frac{(-\ad(\mu_n))^i ( \pi_\lambda A \pi_{\lambda'}
)} {(k-\lambda+\lambda')^{i+1}} ,
\end{align*}
where of course, the sum is finite since $\mu_n$ is nilpotent.

Let $\Theta(z)$ and $\tilde\Theta(z)=\Theta(z)\rho(z)$
be a pair of fundamental solutions, where
$$
\rho(z) = \Id + \sum_{k=1}^\infty z^k \rho_k
$$
is a flat endomorphism of $\CH=T^*M[z,z^{-1}]$. Let $\p_z+z^{-1}R(z)$ and
$\p_z+z^{-1}\tilde{R}(z)$ be the associated endomorphisms of $\CH$. Then we
have
$$
\tilde{R}_k = (k-\ad(\mu)) G_k + \sum_{\ell=1}^{k-1} \bigl( G_\ell
R_{k-\ell} - \tilde{R}_{k-\ell} G_\ell \bigr) + R_k .
$$
Setting
$$
G_k = - \sum_{\lambda-\lambda'\ne k} (k-\ad(\mu))^{-1} \biggl(
\sum_{\ell=1}^{k-1} \pi_\lambda \bigl( G_\ell R_{k-\ell} -
\tilde{R}_{k-\ell} G_\ell \bigr) \pi_{\lambda'} + \pi_\lambda R_k
\pi_{\lambda'}  \biggr) ,
$$
we obtain a fundamental solution $\tilde{\Theta}(z)$ such that
$[\mu_s,\tilde{R}_k]=k\tilde{R}_k$.
\end{proof}

\begin{definition}
A \textbf{conformal} fundamental solution is a fundamental solution
satisfying the condition $[\mu_s,R_k]=kR_k$
\end{definition}

Extracting the constant term of the equation
$\Theta(z)\*\delta_z-\tdelta_z\*\Theta(z)=0$, we see that
\begin{equation} \label{U}
\U = [\mu,\M] + \M + R_1 .
\end{equation}
Applying $\p$ to this equation, we see that
\begin{equation} \label{X}
\p\U = [\mu,\X] + \X .
\end{equation}

\subsection{Vector fields on the jet-space of a conformal Dubrovin manifold}
Let $\L_k$ be the vector field on the jet-space $J^\infty_0M$ of a
conformal Dubrovin manifold $(M,\nabla^z,e,E,r)$ defined by the
formula
$$
\L_k = \Res_{z=0} \bigl\<\G(z),z\delta_z^{k+1}\D(z)\bigr\> .
$$
For example, $\L_{-1}=-e$. By the work of Eguchi, Hori and Xiong
\cite{EHX} and Dubrovin and Zhang \cite{DZ}, these vector fields play
a fundamental role in the theory of Gromov-Witten invariants. In this
section, we give simple intrinsic formulas for them.
\begin{theorem}
  Restricted to $\CO_0 \subset \CO_\infty$, the vector field $\L_k$
  equals $-E^{\circ k+1}$. We have
$$
  [\p,\L_k] = \Res_{z=0}\bigl\<e,\delta_z^{k+1}\D(z)\bigr\> ,
$$
and
\begin{equation} \label{Lk}
\L_k = - \sum_{\ell=0}^\infty \biggl( \p^\ell \<e,\U^{k+1}\> +
\sum_{j=1}^\ell \p^{j-1} \Res_{z=0} \bigl\< e ,
\delta_z^{k+1}(\p+z\X)^{\ell-j}\X \bigr\> \biggr) \p_\ell .
\end{equation}
\end{theorem}
\begin{proof}
  By Theorem \ref{G}, $\G_\ell$ vanishes for $\ell\le0$. Since
  $\D(z)u^a=\X du^a$, it follows that
\begin{align*}
  \L_ku^a &= \Res_{z=0} \bigl\< \G(z) , z\delta_z^{k+1}\X du^a \bigr\>
  = \<\G_1,\U^{k+1}\X du^a\> \\
  &= \< \X^*\G_1 , \U^{k+1} du^a \> = - \<e,\U^{k+1}du^a\> = -
  E^{\circ k+1}(u^a) .
\end{align*}
Here, we have used that $[\U,\X]=0$.

Choose a conformal fundamental solution $\Theta(z)$ on $M$. We may
replace $\G(z)$ and $\D(z)$ in the formula for $\L_k$ by
$\Theta^*(z)\g(z)$ and $\Theta^{-1}(z)\d(z)$, and rewrite $\L_k$ as
$$
\L_k = \Res_{z=0}
\bigl\<\Theta^*(z)\g(z),z\delta_z^{k+1}\Theta^{-1}(z)\d(z)\bigr\> =
\Res_{z=0} \bigl\<\g(z),z\tdelta_z^{k+1}\d(z)\bigr\> .
$$
From this formula, and the fact that $[\p,\tdelta_z]=0$ and
$[\p,\d(z)]=0$, we see that
\begin{align*}
[\p,\L_k] &= \Res_{z=0} \bigl\<\p\g(z),z\tdelta_z^{k+1}\d(z)\bigr\>
= \Res_{z=0} \bigl\< \Theta^*(z)^{-1}e , \tdelta_z^{k+1}\d(z)\bigr\> \\
&= \Res_{z=0} \bigl\< e , \delta_z^{k+1}\D(z)\bigr\> .
\end{align*}

If $X$ is a vector field on the jet space $J^\infty M$, we have
$$
X(u^a_\ell) = \p^\ell(Xu^a) - \sum_{j=1}^\ell \p^{j-1} \cdot [\p,X]
u^a_{\ell-j} .
$$
In the special case $X=\L_k$, this gives
$$
\L_k u^a_\ell = - \p^\ell(\L_ku^a) - \sum_{j=1}^\ell \p^{j-1} \Res_{z=0}
\< e , \delta_z^{k+1}\D(z) u^a_{\ell-j} \> ,
$$
from which \eqref{Lk} follows.
\end{proof}

Using the noncommutative Fa\`a di Bruno polynomials $\PP_n(z)$, we may
rewrite \eqref{Lk} in the elegant form
$$
\L_k = - \sum_{\ell=0}^\infty \sum_{j=0}^\ell \p^j \Res_{z=0} \bigl\<
e , \delta_z^{k+1} \PP_{\ell-j}(z) \bigr\> \p_\ell .
$$
The following corollary is an immediate consequence.
\begin{corollary}
  The vector field $\L_k$ preserves the filtration of $\CO_\infty$ by
  subalgebras $\CO_j$.
\end{corollary}

As an illustration of the utility of the formula \eqref{Lk}, let us
calculate $\L_k\X$; we will need this calculation later in the
discussion of Virasoro constraints in genus $1$.
\begin{proposition} \label{EkX}
$$
\L_k\X = \sum_{\ell=0}^k [\X,\U^\ell\mu\U^{k-\ell}] - \sum_{\ell=0}^k
\bigl( \U^\ell\mu\U^{k-\ell} \bigr) {}_e^a \A{a} \X - \tfrac{3}{2} (k+1)
\X\U^k
$$
\end{proposition}
\begin{proof}
We have $\X=u^b_1\A{b}$, hence $\p_{0,a}\X=\p\A{a}$, $\p_{1,a}\X=\A{a}$,
while $\p_{n,a}\X=0$ for $n>1$. It follows that
\begin{align*}
\L_k\X &= - \<e,\U^{k+1}\> \, \p_0\X - \p\<e,\U^{k+1}\> \, \p_1\X -
\Res_{z=0} \bigl\< e , \delta_z^{k+1}\X \bigr\> \, \p_1\X \\
&= - \p\bigl( (\U^{k+1}){}_e^a \A{a} \bigr) - \Res_{z=0} \bigl(
\delta_z^{k+1}\X \bigr){}_e^a \A{a} \\
&= - \p(\U^{k+1}) - \Res_{z=0} \bigl( \delta_z^{k+1} \bigr){}_e^a \A{a}\X ,
\end{align*}
since $[\A{X},\A{Y}]=0$. By \eqref{X}, we have
$$
\p(\U^{k+1}) = \sum_{\ell=0}^k \U^\ell ( \X + [ \mu , \X ] ) \U^{k-\ell}
= - (k+1) \X\U^k + \sum_{\ell=0}^k [ \X , \U^\ell\mu\U^{k-\ell} ] .
$$
We also have
\begin{align*}
\Res_{z=0} \bigl( \delta_z^{k+1} \bigr){}_e^a \A{a}\X &= \sum_{\ell=0}^k
\bigl( \U^\ell \bigl( \mu + \half \bigr) \U^{k-\ell} \bigr){}_e^a \A{a}\X \\
&= \half (k+1) \bigl( \U^k \bigr){}_e^a \A{a}\X + \sum_{\ell=0}^k \bigl(
\U^\ell\mu\U^{k-\ell} \bigr){}_e^a \A{a}\X .
\end{align*}
The proposition follows, since $\bigl( \U^k \bigr){}_e^a\A{a}=\U^k$.
\end{proof}

\begin{corollary}
In terms of the socle one-form $\om$ of \eqref{socle},
$$
\L_k\log\det(\X) = - \sum_{\ell=0}^k \bigl< e ,
\U^\ell(\mu+\tfrac{3}{2})\U^{k-\ell} \om \bigr\> .
$$
\end{corollary}
\begin{proof}
We see from Proposition \ref{EkX} that
\begin{multline*}
\L_k\log\det(\X) = \Tr(\X^{-1}\L_k\X) \\
= \sum_{\ell=0}^k \Tr[\X,\X^{-1}\U^\ell\mu\U^{k-\ell}] - \sum_{\ell=0}^k
\bigl( \U^\ell\mu\U^{k-\ell} \bigr) {}_e^a \Tr(\A{a}) - \tfrac{3}{2} (k+1)
\Tr(\U^k) .
\end{multline*}
The corollary follows, since $\Tr[\X,\X^{-1}\U^\ell\mu\U^{k-\ell}]=0$ and
$\Tr(\U^k)=\<e,\U^k\om\>$. 
\end{proof}

The vector fields $\L_k$ form a Lie algebra isomorphic to the Lie
algebra of vector fields on the line, under the correspondence which
associates to $\L_k$ the vector field $-z\p_z^{k+1}$. These relations
are sometimes called the Virasoro relations.
\begin{theorem} \label{virasoro}
$[\L_j,\L_k] = (j-k) \L_{j+k}$
\end{theorem}
\begin{proof}
We have
\begin{multline*}
\Bigl[ \bigl\<\g(y),y\tdelta_y^{j+1}\d(y)\bigr\> ,
\bigl\<\g(z),z\tdelta_z^{k+1}\d(z)\bigr\> \Bigr] \\
\begin{aligned}
&= \biggl\< \g(y) , y\tdelta_y^{j+1} \*
\frac{\Theta(y)\Theta^{-1}(z)}{z-y} \* z\tdelta_z^{k+1} \d(z) \biggr\> \\
&- \biggl\< \g(z) , z\tdelta_z^{k+1} \* \frac{\Theta(z)\Theta^{-1}(y)}{y-z}
\* y\tdelta_y^{j+1} \d(z) \biggr\> \\
&= \biggl\< \G(y) , y\delta_y^{j+1} \* \frac{1}{z-y} \* z\delta_z^{k+1} \D(z)
\biggr\>
- \biggl\< \G(z) , z\delta_z^{k+1} \* \frac{1}{y-z} \* y\delta_y^{j+1} \D(z)
\biggr\> \\
&= \frac{1}{z-y} \bigl\< (\delta_y^*)^{j+1}y\G(y) , z\delta_z^{k+1} \D(z)
\bigr\> - \frac{1}{y-z} \bigl\< (\delta_z^*)^{k+1} z\G(z)
, y\delta_y^{j+1} \D(z) \bigr\> .
\end{aligned}
\end{multline*}
The coefficients of $y^\ell$ and $z^\ell$ in $(\delta_y^*)^{j+1}y\G(y)$ and
$(\delta_z^*)^{k+1}z\G(z)$ vanish for $\ell\ge0$; it follows that
\begin{align*}
\Res_{y=0}\frac{(\delta_y^*)^{j+1}y\G(y)}{z-y} &= (\delta_z^*)^{j+1}z\G(z) , &
\Res_{z=0} \frac{(\delta_z^*)^{k+1} z\G(z)}{y-z} &= (\delta_y^*)^{k+1} y\G(y) ,
\end{align*}
and hence that
\begin{multline*}
[\L_j,\L_k] = \Res_{y=0} \Res_{z=0} \Bigl[
\bigl\<\g(y),y\tdelta_y^{j+1}\d(y)\bigr\> ,
\bigl\<\g(z),z\tdelta_z^{k+1}\d(z)\bigr\> \Bigr] \\
\begin{aligned}
&= \Res_{z=0} \bigl\< (\delta_z^*)^{j+1}z\G(z) , z\delta_z^{k+1}\D(z)
\bigr\> - \Res_{y=0} \bigl\< (\delta_y^*)^{k+1} y\G(z) ,
y\delta_y^{j+1}\D(y) \bigr\> \\
&= \Res_{z=0} \bigl\< \G(z) , \bigl(z\delta_z^{j+1}\*z\delta_z^{k+1}
- z\delta_z^{k+1}\*z\delta_z^{j+1} \bigr) \D(y) \bigr\> .
\end{aligned}
\end{multline*}
The result follows, since $z\delta_z^{j+1}\*z\delta_z^{k+1} =
z^2\delta_z^{j+k+2}+(j+1)z\delta_z^{j+k+1}$.
\end{proof}

As a corollary, we see that
$$
[E^{\circ j},E^{\circ k}] = (k-j) E^{\circ(j+k-1)} .
$$
This relation was first conjectured for Frobenius manifolds by
Dubrovin and Zhang \cite{DZ}, and proved by Hertling and Manin
\cite{HM}, using their theory of $F$-manifolds.

\subsection*{Relationship to prior results}
The theory of conformal Dubrovin connections is due to Dubrovin. The
study of the Virasoro vector fields on the jet-space of a Frobenius
manifold was initiated by Dubrovin and Zhang \cite{DZ}; the
presentation in Section~3.2 is based on the methods of Section~5 of
\cite{virasoro} (for the action of these vector fields on functions of
the underlying Dubrovin manifold) and of \cite{EGX} (for the
commutator with the vector field $\p$).

\section{Topological recursion relations for Gromov-Witten invariants and
the jet-space} \label{definition}

The genus $0$ Gromov-Witten invariants of a smooth projective variety
$X$ define a conformal Dubrovin manifold with conformal fundamental
solution, called the small phase space of $X$. In this section, we
recall the construction of the small phase space, following Dubrovin
\cite{D}. We then show the way in which the theory of topological
recursion relations (Eguchi and Xiong \cite{EX}, Kontsevich and Manin
\cite{KM}, and our own work \cite{genus2}) takes a simpler form when
rewritten using the intrinsic geometry of the jet-space of this
Dubrovin manifold; this is the main new result of this paper.

\subsection{Stable maps and Gromov-Witten invariants}
Let us recall the definition of the Gromov-Witten invariants of a
projective manifold over $\C$; see Cox and Katz \cite{CK} and Manin
\cite{M} for more detailed expositions. The definition which we
outline is the one which works in the setting of algebraic geometry:
Gromov-Witten invariants have also been defined for compact symplectic
manifolds, using entirely different techniques. The Gromov-Witten
invariants of a projective manifold reflect the intersection theory of
Kontsevich's moduli spaces of stable maps $\Mbar_{g,n}(X,\beta)$,
whose definition we now recall.

Let $X$ be a projective manifold of dimension $r$. A \textbf{prestable map}
$$\textstyle
(f:C\to X,z_1,\dots,z_n)
$$
of genus $g\ge0$ and degree $\beta\in H_2(X,\Z)$ with $n$ marked points
consists of the following data:
\begin{enumerate}
\item a connected projective curve $C$ of arithmetic genus
$g=h^1(C,\CO_C)$, whose only singularities are ordinary double points,
\item $n$ distinct smooth points $(z_1,\dots,z_n)$ of $C$;
\item an algebraic map $f:C\to V$, such that the degree of $f$, that
is, the cycle $f_*[C]\in H^2(X,\Z)$, equals $\beta$.
\end{enumerate}
If $\tilde{C}$ is the normalization of $C$, the special points in
$\tilde{C}$ are the inverse images of the singular and marked points of
$C$. (Note that the degree of $f:C\to X$ equals $0$ if and only if its
image is a single point.)

A prestable map $(f:C\to X,z_1,\dots,z_n)$ is \textbf{stable} if it has no
infinitesimal automorphisms fixing the marked points. The condition of
stability is equivalent to the following: each irreducible component of
$\tilde{C}$ of genus $0$ on which $f$ has degree $0$ has at least $3$
special points, while each irreducible component of $\tilde{C}$ of genus
$1$ on which $f$ has degree $0$ has at least $1$ special point. In
particular, there are no stable maps of genus $g$ and degree $0$ with $n$
marked points unless $2(g-1)+n>0$.

The moduli stack of $n$-pointed stable maps $\Mbar_{g,n}(X,\beta)$ is the
classifying stack for stable maps of genus $g$ and of degree $\beta$; it is
a complete Deligne-Mumford stack, though not in general smooth (Behrend and
Manin \cite{BM}). The definition of Gromov-Witten invariants is based on
the study of $\Mbar_{g,n}(X,\beta)$. Let $\ev_i:\Mbar_{g,n}(X,\beta)\to X$,
$1\le i\le n$, be evaluation at the $i$th marked point:
$$
\ev_i : \Mbar_{g,n}(X,\beta) \ni { \textstyle (f:C\to X,z_1,\dots,z_n) }
\mapsto f(z_i) \in X .
$$

Let $(f:C\to X,z_1,\dots,z_{n+N})$ be a stable map, representing a
point of $\Mbar_{g,n+N}(X,\beta)$. Forgetting the last $N$ points, we
obtain a map $(f:C\to X,z_1,\dots,z_n)$ which may not be stable, owing
to the existence of rational components of $C$ on which the sheaf
$f^*TX\o\CO_C(-(z_1+\dots+z_n))$ has non-vanishing sections. On such a
component, the map $f$ is constant; thus, there is a stable map
$(f:C\to X,z_1,\dots,z_n)^\stab$ obtained by contracting each of these
components of $C$. There is a morphism, constructed by Behrend and
Manin \cite{BM}, which sends $(f:C\to X , z_1,\dots,z_{n+N})$ to
$(f:C\to X , z_1,\dots,z_n)^\stab$, and which we denote
$$
\pi_{n,N} : \Mbar_{g,n+N}(X,\beta) \to \Mbar_{g,n}(X,\beta) .
$$

In the case $N=1$, this construction yields the \textbf{universal curve}
\begin{align} \label{universal}
\pi = \pi_{n,1} : \Mbar_{g,n+1}(X,\beta) \to \Mbar_{g,n}(X,\beta) .
\end{align}
The fibre of $\pi$ at a stable map $(f:C\to X,z_1,\dots,z_n)$ is the curve
$C$; $f=\ev_{n+1}:\Mbar_{g,n+1}(X,\beta) \to X$ is the universal stable map.

The sheaf $R^1\pi_*f^*TX$ on $\Mbar_{g,n}(X,\beta)$ is called the
obstruction sheaf of $\Mbar_{g,n}(X,\beta)$. If it vanishes, the
Grothendieck-Riemann-Roch theorem implies that the stack
$\Mbar_{g,n}(X,\beta)$ is smooth, of dimension
\begin{equation} \label{vdim}
\vdim\Mbar_{g,n}(X,\beta) = (3-r)(g-1) + \tint_\beta c_1(X) + n ;
\end{equation}
this is called the \textbf{virtual dimension} of $\Mbar_{g,n}(X,\beta)$.
This hypothesis is rarely true; however, there is an algebraic cycle
$$
[\Mbar_{g,n}(X,\beta)]^\virt \in
H_{2\vdim\Mbar_{g,n}(X,\beta)}(\Mbar_{g,n}(X,\beta),\Q) ,
$$
the \textbf{virtual fundamental class}, which stands in for
$[\Mbar_{g,n}(X,\beta)]$ in the general case. The most important
property of the virtual fundamental class is the formula
\begin{equation} \label{stabvirt}
[\Mbar_{g,n+1}(X,\beta)]^\virt = \pi^![\Mbar_{g,n}(X,\beta)]^\virt
\end{equation}
for integration over the fibres of the stabilization map.

If the obstruction sheaf $R^1\pi_*f^*TX$ is locally free of rank $e$, the
moduli stack $\Mbar_{g,n}(X,\beta)$ is smooth of dimension
$\vdim\Mbar_{g,n}(X,\beta)+e$, and the virtual fundamental class satisfies
the equation
$$
[\Mbar_{g,n}(X,\beta)]^\virt = c_e(R^1\pi_*f^*TX) \cap
[\Mbar_{g,n}(X,\beta)] .
$$

The universal curve \eqref{universal} over $\Mbar_{g,n}(X,\beta)$ has $n$
canonical sections
$$
\sigma_i : \Mbar_{g,n}(X,\beta)\to\Mbar_{g,n+1}(X,\beta) ,
$$
corresponding to the $n$ marked points of the curve $C$. Consider the line
bundles
$$
\Om_i = \sigma_i^*\om , \quad 1\le i\le n ,
$$
where 
$\om = \om_{\Mbar_{g,n+1}(X,\beta)/\Mbar_{g,n}(X,\beta)}$
is the relative dualizing sheaf of the universal curve. Thus, the fibre of
$\Om_i$ at the stable map $f:C\to X$ equals the cotangent line $T^*_{z_i}C$
of $C$ at the $i$th marked point $z_i$ of $C$. Let
$$
\Psi_i=c_1(\Om_i) \in H^2(\Mbar_{g,n}(X,\beta),\Z)
$$
be the Chern class of $\Om_i$.

The Gromov-Witten invariants of a projective manifold are defined by
integrating tautological cohomology classes against the virtual fundamental
class $[\Mbar_{g,n}(X,\beta)]^\virt$. Given rational cohomology classes
$x_1,\dots,x_n$ on $X$, we define the Gromov-Witten invariants by the
formula
$$
\< \tau_{k_1}(x_1) \dots \tau_{k_n}(x_n) \>^X_{g,\beta} =
\int_{[\Mbar_{g,n}(X,\beta)]^\virt} \Psi_1^{k_1} \dots \Psi_n^{k_n}
\ev_1^*x_1 \dots \ev_n^*x_n \in \Q .
$$
Perhaps the most important property of the Gromov-Witten invariants is that
they are invariant under deformation of the variety $X$.

The \textbf{Novikov ring} $\Nov$ of $X$ is the commutative graded algebra
consisting of all formal sums
$$
\sum_{\beta\in H_2(X,\Z)} a_\beta q^\beta
$$
such that for all $C>0$, the set of $\beta\in H_2(X,\Z)$ such that
$a_\beta\ne0$ and $\tint_\beta\om<C$ is finite; the product is defined
by $q^{\beta_1}q^{\beta_2}=q^{\beta_1+\beta_2}$, and
the grading is defined by $|q^\beta|=-2c_1(X)\cap\beta$. For example,
for $\CP^d$, we have $\Nov=\Q[q]$, where $|q|=-2(d+1)$. By working
over the Novikov ring, we may combine the Gromov-Witten invariants in
different degrees into a single generating function:
$$
\<\tau_{k_1}(x_1)\dots\tau_{k_n}(x_n)\>^X_g = \sum_{\beta\in H_2(X,\Z)}
q^\beta \<\tau_{k_1}(x_1)\dots\tau_{k_n}(x_n)\>^X_{g,\beta} .
$$

\subsection{The small and large phase spaces} Let
$$
\{ \gamma_a \in H^{p_a,q_a}(X) \mid a\in A\}
$$
be a homogeneous basis of the Dolbeault cohomology of $X$, such that
for a distinguished element $e\in A$, $\gamma_e=1$. The small phase space
$\H$ of $X$ is the formal neighbourhood of $0$ in the vector space
$H^\bull(X,\C)$. Denote by $u^a$ the coordinates on $\H$ dual to the basis
$\gamma_a$. Let $\eta$ be the flat metric on $\H$ associated to the
Poincar\'e form on $H^\bull(X,\C)$, with components
$$
\eta_{ab} = \tint_X \gamma_a \cup \gamma_b ,
$$
and by $\eta^{ab}$ the inverse matrix of $\eta_{ab}$. 

Dubrovin \cite{D} shows that $\H$ is a conformal Dubrovin manifold, and
construct from the genus $0$ Gromov-Witten invariants of $X$ a conformal
fundamental solution $\Theta(z)$ on $\H$. The flat connection $\nabla$ of
this Dubrovin manifold is the Levi-Civit\`a connection associated to the
flat metric $\eta$. We defer the definition of the remaining geometric
structures on $\H$, namely the tensor $\AA$, the identity and Euler vector
fields $e$ and $E$ and the fundamental solution $\Theta$, to below.

The large phase space $\H_\infty$ of $X$ is the formal neighbourhood of $0$
in the vector space $H^\bull(X\times\CP^\infty,\C)$; it is a formal
manifold with coordinates $\{t^a_k \mid a\in A,k\ge0\}$, where $t^a_k$ has
degree $-p_a-2k$. (If $X$ has cohomology of odd degree, $\H$ and
$\H_\infty$ are actually supermanifolds; this detail does not materially
change the theory.)

The genus $g$ potential $\PHI_g$ of $X$ is the function on the large phase
space $\H_\infty$ given by the formula
\begin{equation} \label{Phi}
\PHI^X_g = \sum_{n=0}^\infty \frac{1}{n!} \sum_{\substack{ k_1\dots k_n \\
a_1 \dots a_n }} t^{a_1}_{k_1} \dots t^{a_n}_{k_n} \< \tau_{k_1,a_1} \dots
\tau_{k_n,a_n} \>^X_g ,
\end{equation}
where we write $\tau_{k,a}$ in place of $\tau_k(\gamma_a)$.

Denote the constant coefficient vector field $\p/\p t^a_k$ on $\H_\infty$
by $\d_{k,a}$, and $\d_{0,e}$ by $\p$. The partial derivatives of the
potential $\PHI_g$ are denoted
$$
\<\< \tau_{k_1,a_1} \dots \tau_{k_n,a_n} \>\>^X_g = \d_{k_1,a_1} \dots
\d_{k_n,a_n} \PHI^X_g .
$$
In particular, $\<\tau_{k_1,a_1} \dots \tau_{k_n,a_n}
\>^X_g$ is the value at $0\in\H_\infty$ of the function $\<\<
\tau_{k_1,a_1} \dots \tau_{k_n,a_n} \>\>^X_g$. Explicitly, we have
\begin{multline*}
\<\<\tau_{k_1,a_1}\dots\tau_{k_n,a_n}\>\>^X_g = \sum_{\beta\in H_2(X,\Z)}
q^\beta \sum_{N=0}^\infty \frac{1}{N!} \sum_{\substack{b_1\dots b_N \\
k_{n+1}\dots k_{n+N} }} t^{b_1}_{k_{n+1}} \dots t^{b_N}_{k_{n+N}} \\
\int_{[\Mbar_{g,n+N}(X,\beta)]^\virt} \Psi_1^{k_1}\dots\Psi_{n+N}^{\ell_N}
\ev_1^*\gamma_{a_1}\dots\ev_n^*\gamma_{a_n}\ev_{n+1}^*\gamma_{b_1}
\dots\ev_{n+N}^*\gamma_{b_N} .
\end{multline*}

\subsection{Dual graphs and their associated cycles}
The moduli space of stable maps $\Mbar_{g,n}(X,\beta)$ is stratified
by boundary strata, which are most conveniently parametrized by dual
graphs.

The dual graph $C(G)$ of a pointed prestable map $(f:C\to X,z_1,\dots,z_n)$
is a graph $G=G(C)$ with one vertex $v$ for each component $C(v)$ of the
normalization $\tilde{C}$ of the curve $C$, labelled by the genus $g(v)$ of
this component and the degree $\beta(v)$ of the restriction $f(v):C(v)\to
X$ of $f$ to $C(v)$.

The edges of the dual graph $G(C)$ correspond to double points of the curve
$C$; the two ends of an edge are attached to the vertices associated to the
components on which the two branches of the double point lie. (If both
branches lie in the same component of $C$, then the edge is a loop.)

Finally, to each marked point $z_i$ of the curve corresponds a leg of the
graph, labelled by $i$, at the vertex corresponding to the component of $C$
on which $z_i$ lies (which is uniquely determined, since $z_i$ is a smooth
point).

In drawing dual graphs, we denote vertices of genus $0$ either by a solid
circle $\begin{picture}(10,8)(0,0)\put(5,3){\circle*{4}}\end{picture}$ or
leave them unmarked and vertices of genus $g>0$ by
$\begin{picture}(10,8)(0,0)\put(5,3){\circle{7}}
\put(3,2){$\scriptscriptstyle g$}\end{picture}$.

Define the genus $g(G)$ of a dual graph $G$ to be the sum of the
genera $g(v)$ over the vertices $v\in\Vert(G)$ of $G$ and the first
Betti number of the graph $G$. Then $g(G)$ equals the arithmetic genus
of the curve $C$. Define the degree $\beta(G)$ of a dual graph $G$ to
be the sum of the degrees $\beta(v)$ over the vertices $v$ of $G$. The
valence $n(G)$ equals the number of external legs of $G$.

A dual graph is \textbf{stable} if each vertex $v$ such that $\beta(v)=0$
satisfies the additional condition $2(g(v)-1)+n(v)>0$; a pointed prestable
curve is stable if and only if the associated dual graph is stable. The set
of isomorphism classes of stable dual graphs $G$ of fixed genus $g(G)$,
valence $n(G)$, and degree $\beta(G)$ is finite.

If $G$ is a stable graph of genus $g$, valence $n$ and degree $\beta$,
let $\M(G)\subset\Mbar_{g,n}(X,\beta)$ be the moduli stack of stable
maps with dual graph $G$, and let $\Mbar(G)$ be its closure in
$\Mbar_{g,n}(X,\beta)$. Let $\Mtilde(G)$ be the product
$$
\Mtilde(G) = \prod_{v\in\Vert(G)} \Mbar_{g(v),n(v)}(X,\beta(v)) ,
$$
with virtual fundamental class
$$
[\Mtilde(G)]^\virt = \prod_{v\in\Vert(G)}
[\Mbar_{g(v),n(v)}(X,\beta(v))]^\virt .
$$
There is a natural covering map
$$
\pi(G) : \Mtilde(G) \longrightarrow \Mbar(G) ,
$$
with covering group $\Aut(G)$. Let $[G]=[\Mbar(G)]^\virt$ be the cycle
$$
[\Mbar(G)]^\virt = \frac{\pi(G)_*[\Mtilde(G)]^\virt}{|\Aut(G)|} .
$$
In formulas, we will often symbolize the cycle $[G]$ by the dual graph
$G$ itself.

A dual graph $G$ with one edge determines a Cartier divisor $D(G)$
supported by $\Mbar(G)$. Axiom V for virtual fundamental classes in
Behrend \cite{B} implies that
$$
D(G) \cap [\Mbar_{g(G),n(G)}(X,\beta(G))]^\virt = [G] .
$$

The following formula is due to Witten \cite{Witten}.
\begin{proposition} \label{psi-pull-back-V}
  Let $\pi:\Mbar_{g,n+1}(X,\beta) \to \Mbar_{g,n}(X,\beta)$ be the
  universal curve \eqref{universal}. For $1\le i\le n$, let
  $D_i=D(G_i)$ be the Cartier divisor on $\Mbar_{g,n+1}(X,\beta)$
  associated to the dual graph
$$
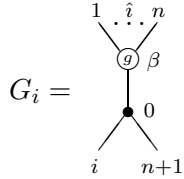
\begin{picture}(90,65)(65,730)
\put(90,760){$G_i = $}
\put(135,771){\line( 0,-1){ 15}}
\put(135,755){\circle*{4}}
\put(141,753){$\scriptstyle 0$}
\put(134,754){\line(-3,-4){ 10}}
\put(136,754){\line( 3,-4){ 10}}
\put(129,788){$\dots$}
\put(121,791){$\scriptstyle 1$}
\put(134,791){$\scriptstyle\hat\imath$}
\put(144,791){$\scriptstyle n$}
\put(135,775){\circle{8}}
\put(133,774){$\scriptscriptstyle g$}
\put(142,772){$\scriptstyle \beta$}
\put(132,778){\line(-3, 4){ 8}}
\put(138,778){\line( 3, 4){ 8}}
\put(121,732){$\scriptstyle i$}
\put(140,732){$\scriptstyle n+1$}
\end{picture}
$$
Then on $\Mbar_{g,n+1}(X,\beta)$, we have the formula $\pi^*\Psi_i =
\Psi_i + D_i$.
\end{proposition}

\subsection{The puncture, dilaton and divisor equations on the large phase
  space} For $\om\in H^2(X,\C)$, let $R_a^b(\om)$ be the matrix of
multiplication by $\om$ on $H(X)$:
$$
\om \cup \gamma_a = R_a^b(\om) \gamma_b .
$$
Combining Proposition \ref{psi-pull-back-V} with \eqref{stabvirt},
we obtain the \textbf{puncture}, \textbf{divisor} and \textbf{dilaton}
equations:
\begin{align}
\<\tau_{0,e}\tau_{k_1,a_1}\dots\tau_{k_n,a_n}\>_{g,\beta}^X &= \sum_{i=1}^n
\<\tau_{k_1,a_1}\dots\tau_{k_i-1,a_i} \dots \tau_{k_n,a_n}\>_{g,\beta}^X ,
\label{puncture} \\
\<\tau_0(\om)\tau_{k_1,a_1}\dots\tau_{k_n,a_n}\>_{g,\beta}^X &=
\sum_{i=1}^n R_{a_i}^b(\om) \<\tau_{k_1,a_1}\dots\tau_{k_i-1,b} \dots
\tau_{k_n,a_n}\>_{g,\beta}^X \label{divisor} \\
&\quad + \tint_\beta\om \*
\<\tau_{k_1,a_1}\dots\tau_{k_n,a_n}\>_{g,\beta}^X , \notag \\
\<\tau_{1,e}\tau_{k_1,a_1}\dots\tau_{k_n,a_n}\>_{g,\beta}^X &= (2g-2+n)
\<\tau_{k_1,a_1}\dots\tau_{k_n,a_n}\>_{g,\beta}^X .
\label{dilaton}
\end{align}
Owing to the nonexistence of the stabilization maps
$\pi:\Mbar_{0,3}(X,0)\to\Mbar_{0,2}(X,0)$ and
$\pi:\Mbar_{1,1}(X,0)\to\Mbar_{1,0}(X,0)$, we have the following
exceptional cases:
\begin{align*}
\<\tau_{0,e}\tau_{0,a}\tau_{0,b}\>_{0,0}^X &= \eta_{ab} , &
\<\tau_{1,e}\>_{1,0}^X &= \tfrac{1}{24} \chi(X) , \\
\<\tau_0(\om)\tau_{0,a}\tau_{0,b}\>_{0,0}^X &= R_{ab}(\om) , &
\<\tau_0(\om)\>_{1,0}^X &= \tfrac{1}{24} \tint_X \om \cup c_{r-1}(X) .
\end{align*}

Introduce the vector fields
\begin{align*}
e &= - \sum_{k=0}^\infty \g_{k+1}^a \p_{k,a} =
\p - \sum_{k=0}^\infty t_{k+1}^a \p_{k,a} , &
\Dil &= - \sum_{k=0}^\infty \g_k^a \p_{k,a}
= \p_{1,e} - \sum_{k=0}^\infty t_k^a \p_{k,a}
\end{align*}
on the large phase space $\H_\infty$. The puncture and dilaton equations
are equivalent to the differential equations
\begin{align*}
e\PHI_g^X &= \begin{cases} \half\eta_{ab} t_0^a t_0^b , & g=0 , \\ 0 , & g>0 ,
\end{cases} &
\Dil\PHI_g^X &= \begin{cases} \tfrac{1}{24} \chi(X) , & g=1 , \\
(2g-2)\PHI_g^X , & g\ne1 ,
\end{cases}
\end{align*}
for the Gromov-Witten potentials $\PHI_g^X$.

We may identify the large phase space $\H_\infty$ with the jet space
$J^\infty_0\H$ of the small phase space $\H$. Consider the map $u$ from the
large phase space $\H_\infty$ to the small phase space $\H$ given by the
formula
$$
u^a = \eta^{ab} \, \p\<\< \d_{0,b} \>\>^X_0 .
$$
Applying the vector field $\p_{0,b}$ to the genus $0$ puncture equation
$e\PHI_0^X=\half\eta_{ab} t_0^a t_0^b$, we see that
$$
u^a = t^a_0 + \sum_{k=0}^\infty t^b_{k+1} \d_{k,b} u^a .
$$
It follows that the map $u:\H_\infty\to\H$ is a submersion at $0$;
this submersion has a section, which identifies $\H$ with the
submanifold of $\H_\infty$ along which the coordinates $t^a_k$, $k>0$,
vanish.

The following result shows that the jet coordinates $\p^nu^a$ form a
coordinate system on the large phase space, and that the origin of the
large phase space corresponds to the basepoint of the jet-space
$J^\infty_0\H$.
\begin{proposition}
We have $\p^nu^a(0)=\delta_{n,1}\delta_e^a$, and
$\d_{k,b}(\p^nu^a)(0)=\delta_{n,k}\delta_b^a$ if $k\le n$.
\end{proposition}
\begin{proof}
We argue by induction on $n=0$; we have already seen that the result is
true for $n=0$. Expanding the equation $e^n(u^a)=\delta_{n,1}\delta_e^a$,
which holds for $n>0$, we see that
\begin{equation} \label{pnu}
\p^nu^a = \delta_{n,1}\delta_e^a + \sum_{i=0}^{n-1} (-1)^i \tbinom{n}{i}
\biggl( \sum_{k=0}^\infty t_{k+1}^b \d_{k,b} \biggr)^{n-i} \p^iu^a .
\end{equation}
In particular, evaluating at $0\in\H_\infty$, we see that
$\p^nu^a(0)=\delta_{n,1}\delta_e^a$. Applying the vector field $\d_{k,b}$
to \eqref{pnu}
and evaluating at $0\in\H_\infty$, we see that
$$
\d_{k,b}(\p^nu^a)(0) = \sum_{i=0}^{n-1} (-1)^i \tbinom{n}{i}
\d_{k-n+i,b}(\p^iu^a)(0) .
$$
The result now follows by the induction hypothesis.
\end{proof}

The following theorem is due to Hori \cite{Hori}; we will only need it
in genus $0$, where it yields a construction of an Euler vector field
on the small phase space.
\begin{theorem}
Let $\L_0$ be the vector field
\begin{align*}
\L_0 &= \sum_{k=0}^\infty (k+p_a+\tfrac{1-r}{2}) \g_k^a \d_{k,a} +
\sum_{k=0}^\infty R_a^b \g_{k+1}^a \d_{k,b} \\
&= - \half (3-r) \p_{1,e} +\sum_{k=0}^\infty (k+p_a+\tfrac{1-r}{2})
t_k^a \d_{k,a} - R_0^b \d_{0,b} + \sum_{k=0}^\infty R_a^b t_{k+1}^a
\d_{k,b} ,
\end{align*}
where $R_a^b$ is the matrix $R_a^b(c_1(X))$. We have
$$
0 = \L_0\PHI_g^X + \begin{cases} \half R_{ab}t_0^at_0^b , & g=0 , \\
\frac{1}{48} \int_X \bigl( (3-r) c_r(X) - 2c_1(X) c_{r-1}(X) \bigr) , & g=1
, \\
0 , & g>1 .
\end{cases}$$
\end{theorem}
\begin{proof}
The formula \eqref{vdim} for the dimension of the virtual fundamental class
of $\Mbar_{g,n}(X,\beta)$ implies the following identity:
\begin{multline*}
\sum_{i=1}^n (p_{a_i}+k_i)
\<\tau_{k_1,a_1}\dots\tau_{k_n,a_n}\>_{g,\beta}^X =
\vdim\Mbar_{g,n}(X,\beta) \*
\<\tau_{k_1,a_1}\dots\tau_{k_n,a_n}\>_{g,\beta}^X \\
= \bigl( (3-r)(g-1) + \tint_\beta c_1(X) + n \bigr)
\<\tau_{k_1,a_1}\dots\tau_{k_n,a_n}\>_{g,\beta}^X .
\end{multline*}
To eliminate the dependence on the genus $g$, subtract $\half(3-r)$ times
the dilaton equation \eqref{divisor}; after some rearrangement, this gives
\begin{multline*}
\sum_{i=1}^n \bigl( p_{a_i}+k_i+\tfrac{1-r}{2} \bigr)
\<\tau_{k_1,a_1}\dots\tau_{k_n,a_n}\>_{g,\beta}^X \\ = \half (3-r)
\<\tau_{1,e}\tau_{k_1,a_1}\dots\tau_{k_n,a_n}\>_{g,\beta}^X
+ \tint_\beta c_1(X) \* \<\tau_{k_1,a_1}\dots\tau_{k_n,a_n}\>_{g,\beta}^X .
\end{multline*}
To eliminate the dependence on $\beta$, apply the divisor equation
\eqref{divisor} with $\om=c_1(X)$: this yields
\begin{multline*}
\sum_{i=1}^n \Bigl( \bigl( p_{a_i}+k_i+\tfrac{1-r}{2} \bigr) \,
\<\tau_{k_1,a_1}\dots\tau_{k_n,a_n}\>_{g,\beta}^X 
+ R_{a_i}^b \,
\<\tau_{k_1,a_1}\dots\tau_{k_i-1,b}\dots\tau_{k_n,a_n}\>_{g,\beta}^X \Bigr)
\\
= \half (3-r) \<\tau_{1,e}\tau_{k_1,a_1}\dots\tau_{k_n,a_n}\>_{g,\beta}^X
+ \<\tau_{0,e}\tau_{k_1,a_1}\dots\tau_{k_n,a_n}\>_{g,\beta}^X .
\end{multline*}
Taking into account the exceptional cases in degree $0$, we obtain the
theorem.
\end{proof}

\subsection{Topological recursion relations}
A stable curve is a stable map with target a point; the moduli space
$\Mbar_{g,n}=\Mbar_{g,n}(\pt,0)$ of stable curves of genus $g$ with $n$
marked points is a smooth Deligne-Mumford stack of dimension $3g-3+n$, the
Deligne-Mumford-Knudsen moduli stack of stable curves. Denote the line
bundles $\Om_i$ on $\Mbar_{g,n}=\Mbar_{g,n}(\pt,0)$ by $\om_i$, and their
Chern classes $c_1(\om_i)$ by $\psi_i$.

If $n$ and $N$ are nonnegative integers, let $\Mbar_{g,n,N}(X,\beta)$ be
the moduli space of stable maps $\Mbar_{g,n+N}(X,\beta)$. If $2g-2+n>0$,
the stabilization map
$$
\rho_{n,N} : \Mbar_{g,n,N}(X,\beta) \to \Mbar_{g,n}
$$
is the morphism which takes a stable map $(f:C\to
X,z_1,\dots,z_n,x_1,\dots,x_N)$ to the stabilization
$(C,z_1,\dots,z_n)^\stab$. We abbreviate the pullback $\rho_{n,N}^*\psi_i$
by $\psi_i$.
\begin{proposition} \label{stab}
  For $1\le i\le n$, let $D_i=D(G_i)$ be the Cartier divisor on
  $\Mbar_{g,n,N}(X,\beta)$ associated to the following dual graph:
$$
G_i = 
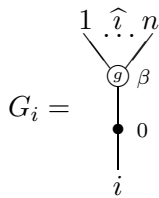
\begin{picture}(30,70)(120,760)
\put(135,771){\line( 0,-1){ 15}}
\put(135,755){\circle*{4}}
\put(142,752){$\scriptstyle 0$}
\put(135,755){\line( 0,-1){ 15}}
\put(133,730){$i$}
\put(135,775){\circle{8}}
\put(133,774){$\scriptscriptstyle g$}
\put(142,772){$\scriptstyle \beta$}
\put(132,778){\line(-3, 4){ 10}}
\put(138,778){\line( 3, 4){ 10}}
\put(129,790){$\dots$}
\put(120,793){$1$}
\put(133,793){$\widehat\imath$}
\put(144,793){$n$}
\end{picture}
$$\vskip25pt
\noindent Then $\Psi_i = \psi_i + D_i$.
\end{proposition}
\begin{proof}
See the proof of Proposition 5 in \cite{genus2}.
\end{proof}

\begin{corollary} \label{Stab}
$$
\frac{1}{1-z\Psi_i} = \frac{1}{1-z\psi_i} \biggl( 1 +
\frac{zD_i}{1-z\Psi_i} \biggr)
$$
\end{corollary}
\begin{proof}
Divide the equation $1 - z\psi_i = 1 - z(\Psi_i-D_i)$ by
$(1-z\psi_i)(1-z\Psi_i)$.
\end{proof}

From Proposition \ref{stab}, we may derive many relations among
Gromov-Witten invariants. For example, using the vanishing of
$\psi_1^{k_1}\psi_2^{k_2}\psi_3^{k_3}$ on the zero-dimensional moduli space
$\Mbar_{0,3}$ when $k_1>0$, we obtain the genus $0$ topological recursion
relation
\begin{equation} \label{trr0}
\<\<\tau_{k_1,a_1}\tau_{k_2,a_2}\tau_{k_3,a_3}\>\>_0 =
\eta^{AB} \<\<\tau_{k_1-1,a_1}\tau_{0,A}\>\>_0
\<\<\tau_{0,B}\tau_{k_2,a_2}\tau_{k_3,a_3}\>\>_0 .
\end{equation}
The following result of Dijkgraaf and Witten \cite{DW}, proved using the
topological recursion relation in genus $0$, may be viewed as an analogue
of Theorem~\ref{Fg} in genus~$0$.
\begin{proposition} \label{DW}
The function $\<\<\tau_{k,a}\tau_{\ell,b}\>\>^X_0$ on the large phase space
is the pullback of a function on the small phase space; that is,
\begin{equation} \label{FPhi}
\<\<\tau_{k,a}\tau_{\ell,b}\>\>^X_0 = \sum_{n=0}^\infty \frac{1}{n!}
u^{a_1} \dots u^{a_n} \< \tau_{k,a} \tau_{\ell,b} \tau_{0,a_1} \dots
\tau_{0,a_n} \>^X_0 .
\end{equation}
\end{proposition}
\begin{proof}
Let $A_{k,a;\ell,b}$ and $B_{k,a;\ell,b}$ be the left and right-hand sides
of \eqref{FPhi}. Since $u^a$ and $t^a_0$ are equal along
$\H\subset\H_\infty$, it follows that $A_{k,a;\ell,b}$ and $B_{k,a;\ell,b}$
are too. We now calculate the derivatives of $A_{k,a;\ell,b}$ and
$B_{k,a;\ell,b}$ with respect to $t^c_m$, $m>0$, using \eqref{trr0}. On the
one hand,
\begin{multline*}
\d_{m,c} B_{k,a;\ell,b} = (\d_{m,c}u^{a_0}) \sum_{n=0}^\infty \frac{1}{n!}
u^{a_1} \dots u^{a_n} \< \tau_{k,a} \tau_{\ell,b}
\tau_{0,a_0} \tau_{0,a_1} \dots \tau_{0,a_n} \>^X_0 \\
\begin{aligned}
&= \eta^{AB} \<\<\tau_{m,c}\tau_{0,e}\tau_{0,A}\>\>^X_0 \sum_{n=0}^\infty
\frac{1}{n!} u^{a_1} \dots u^{a_n} \< \tau_{k,a} \tau_{\ell,b} \tau_{0,B}
\tau_{0,a_1} \dots \tau_{0,a_n} \>^X_0 \\
&= \eta^{AB} \eta^{CD} \<\<\tau_{m-1,c}\tau_{0,C}\>\>^X_0
\,\<\<\tau_{0,D}\tau_{0,e}\tau_{0,A}\>\>^X_0 \\
& \qquad \sum_{n=0}^\infty \frac{1}{n!} u^{a_1} \dots u^{a_n} \< \tau_{k,a}
\tau_{\ell,b} \tau_{0,B} \tau_{0,a_1} \dots \tau_{0,a_n} \>^X_0 \\
&= \eta^{AB} \<\<\tau_{m-1,c}\tau_{0,A}\>\>^X_0 \d_{0,B} B_{k,a;\ell,b} .
\end{aligned}
\end{multline*}
On the other hand, we have
\begin{align*}
\d_{m,c} A_{k,a;\ell,b} &= \<\<\tau_{k,a}\tau_{\ell,b}\tau_{m,c}\>\>^X_0
= \eta^{AB} \<\<\tau_{m-1,c}\tau_{0,A}\>\>^X_0
\<\<\tau_{k,a}\tau_{\ell,b}\tau_{0,B}\>\>^X_0 \\
&= \eta^{AB} \<\<\tau_{m-1,c}\tau_{0,A}\>\>^X_0 \d_{0,B} A_{k,a;\ell,b} .
\end{align*}
Induction in the order of vanishing of
$A_{k,a;\ell,b}-B_{k,a;\ell,b}$ in the variables $\{t^a_m\mid m>0\}$ shows
that the two power series are equal.
\end{proof}

\subsection{The Dubrovin connection on $\H$}
Let $\M$ be the endomorphism of the cotangent bundle of $\H$ with
components $\M_a^b = \eta^{bc} \<\<\tau_{0,a}\tau_{0,c}\>\>^X_0$; by
Proposition \ref{DW}, it is a function on $\H$. We may now define the
Dubrovin connection $\nabla^z=\nabla+z\AA$ on the small phase space
$\H$: the tensor $\AA$ is defined by the formula $\A{X} = X(\M)$.
\begin{proposition}
$\A{X}=\X^{-1}\d_{0,X}\M=(\d_{0,e}\X)^{-1}(\d_{0,X}\M)$
\end{proposition}
\begin{proof}
By \eqref{dp}, we have $\A{a}=(\X^{-1})_a^b\d_{0,b}\M=\X^{-1}\d_{0,a}\M$.
\end{proof}

We now check that $\nabla^z=\nabla+z\AA$ is a Dubrovin connection. The
equation $[\nabla_a,\A{b}]=[\nabla_b,\A{a}]$ follows from the fact
that the tensor $\AA$ is the covariant derivative of $\M$. The
equation $[\A{a},\A{b}]=0$ is equivalent to the
Witten-Dijkgraaf-Verlinde-Verlinde (WDVV) equation
$$
\eta^{AB} \<\<\tau_{0,a}\tau_{0,b}\tau_{0,A}\>\>^X_0
\<\<\tau_{0,B}\tau_{0,c}\tau_{0,d}\>\>^X_0
= \eta^{AB} \<\<\tau_{0,a}\tau_{0,c}\tau_{0,A}\>\>^X_0
\<\<\tau_{0,B}\tau_{0,b}\tau_{0,d}\>\>^X_0 .
$$
This equation is proved by applying the vector field $\tau_{0,b}$ to
the topological recursion relation \eqref{trr0}:
$$
\<\<\tau_{1,a}\tau_{0,c}\tau_{0,d}\>\>_0 =
\eta^{AB} \<\<\tau_{0,a}\tau_{0,A}\>\>_0
\<\<\tau_{0,B}\tau_{0,c}\tau_{0,d}\>\>_0 .
$$
This gives
\begin{multline*}
\eta^{AB} \<\<\tau_{0,a}\tau_{0,b}\tau_{0,A}\>\>_0
\<\<\tau_{0,B}\tau_{0,c}\tau_{0,d}\>\>_0 \\ =
\<\<\tau_{1,a}\tau_{0,b}\tau_{0,c}\tau_{0,d}\>\>_0 -
\eta^{AB} \<\<\tau_{0,a}\tau_{0,A}\>\>_0
\<\<\tau_{0,B}\tau_{0,b}\tau_{0,c}\tau_{0,d}\>\>_0 ;
\end{multline*}
the left-hand side of this equation is symmetric in the indices $b$ and
$c$, since the right-hand side is, and this is precisely the WDVV equation.

The identity vector field $e$ of the Dubrovin connection is the constant
vector field $\p_e$ in the direction of the coordinate $u^e$. The equation
$\A{e}=e\M=\Id$ is obtained by applying the differential operator
$\tau_{0,a}\tau_{0,b}$ to the puncture equation $e\PHI_0^X=\half\eta_{ab}
t_0^a t_0^b=0$.

\subsection{The Euler vector field on $\H$}

Using Hori's equation, we may construct an Euler vector field $E$ on $\H$
which along with the Dubrovin connection $\nabla^z$ and the identity vector
field $e$ makes it into a conformal Dubrovin manifold. To see this, apply
the differential operator $\d_{0,b}\d_{0,c}$ to the equation
$\CL_0\PHI_0^X+\half R_{ab}t_0^at_0^b=0$: we obtain
\begin{equation} \label{hori}
\CL_0\<\<\d_{0,b}\d_{0,c}\>\>_0^X +
(p_b+p_c+1-r)\<\<\d_{0,b}\d_{0,c}\>\>_0^X + R_{bc} = 0 .
\end{equation}
In particular,
$$
\CL_0u^a + (1-p_a)u^a + R_e^a = \CL_0u^a + Eu^a = 0 ,
$$
where $E$ is the vector field
$$
E = \sum_a \bigl( (1-p_a)u^a + R_e^a \bigr) \p_a .
$$
Let $\mu$ be the endomorphism of the cotangent bundle of $\H$ with
components
$$
\mu_a^b = \delta_a^b ( p_a - r/2 ) ,
$$
corresponding to the identification of the constant $r$ in the definition
of the conformal Dubrovin manifold $M$ with the dimension of $X$.
By \eqref{hori}, we see that
$$
\CL_0(\M) + [\mu,\M] + \M + R = 0 .
$$
Applying the vector field $\p$ to \eqref{hori}, we see that
$\CL_0(\X) + [\mu,\X] + \tfrac{3-r}{2} \X = 0$, and hence that
$$
\CL_0(\X^{-1}) + [\mu,\X^{-1}] - \tfrac{3-r}{2} \X^{-1} = 0 .
$$
It follows that
\begin{align*}
E(\A{a}) &= - \CL_0(\X^{-1}\d_{0,a}\M) \\ &= - \CL_0(\X^{-1})\d_{0,a}\M
- \X^{-1}\d_{0,a}\CL_0(\M) - \X^{-1}[\CL_0,\d_{0,a}]\M \\ 
&= \bigl( [\mu,\X^{-1}]\d_{0,a}\M - \tfrac{3-r}{2} \A{a} \bigr) + \bigl(
\X^{-1}\d_{0,a}[\mu,\M] + \A{a} \bigr) + \bigl( p_a + \tfrac{1-r}{2} \bigr)
\A{a} \\
&= [\mu,\A{a}] + p_a \A{a} .
\end{align*}
In other words, $[E,\AA_{ab}^c] = (p_a+p_b-p_c)\AA_{ab}^c$; this implies
\eqref{E}, since on the one hand,
$$
[E,\p_a\circ\p_b] = [E,\AA_{ab}^c\p_c] = \bigl( E(\AA_{ab}^c) - (1-p_c)
\bigr) \p_c ,
$$
while on the other hand,
$$
[E,\p_a]\circ\p_b + \p_a\circ[E,\p_b] + \p_a\circ\p_b
= \bigl( - (1-p_a) - (1-p_b) + 1 \bigr) \AA_{ab}^c\p_c .
$$
This establishes that $E$ is an Euler vector field on the small phase space
$\H$.

Observe that Hori's vector field $\CL_0$ agrees with the vector field
$\CL_0=\Res\<\G(z),z\delta_z\D(z)\>$ associated to the conformal
Dubrovin manifold $(\H,\nabla^z,e,E)$. Historically, Hori's equation
was an important step in the discovery of the Lie algebra of vector
fields $\CL_k$.

\subsection{The fundamental solution on $\H$}
It is remarkable that the small phase space $\H$ has a canonical conformal
fundamental solution; this was first shown by Dubrovin (see \cite{D},
Chapter 6).
\begin{proposition}
The matrix
$$
\Theta_a^b(z) = \delta_a^b + \eta^{bc} \sum_{k=0}^\infty z^{k+1}
\<\<\tau_{k,a}\tau_{0,c}\>\>^X_0
$$
is a conformal fundamental solution on $\H$.
\end{proposition}
\begin{proof}
The equation $\p_a\Theta_k = \Theta_{k-1}\A{a}$ is an instance of the genus
$0$ topological recursion relation:
$$
\<\<\tau_{k,a}\tau_{0,b}\tau_{0,c}\>\>^X_0 = \eta^{AB}
\<\<\tau_{k-1,a}\tau_{0,A}\>\>^X_0
\<\<\tau_{0,B}\tau_{0,b}\tau_{0,c}\>\>^X_0 .
$$
It remains to show that the fundamental solution $\Theta(z)$ is
conformal. Applying the differential operator $\d_{k,a}\d_{0,b}$ to Hori's
equation in genus $0$, we see that
$$
\CL_0\<\<\d_{k,a}\d_{0,b}\>\>_0 + (k+1+p_a+p_b-r)
\<\<\d_{k,a}\d_{0,b}\>\>_0 + \delta_{k,0} R_{ab} = 0 .
$$
Multiplying by $z^{k+1}$ and summing over $k$, we see that
$$
E\Theta(z) = z\p_z\Theta(z) + [\mu,\Theta(z)] + zR .
$$
On the other hand, $E\Theta(z)=z\Theta(z)\U$, where $\U=E\M$. From this
equation, we see that the operator $\tdelta_z$ equals
$$
\tdelta_z = \p_z + z^{-1} ( \mu + \half) + R .
$$
Since $[\mu,R]=R$, we see that $\Theta(z)$ is conformal; this equation
is the statement that multiplication by $c_1(X)$ in the Dolbeault
cohomology raises degree by $(1,1)$.
\end{proof}

\subsection{Intrinsic formulation of topological recursion relations}
We now come to the main result of this paper, an expression for the
topological recursion relations among Gromov-Witten invariants in
terms of the intrinsic vector fields $\sigma_{k,a}$ on the jet-space
$J^\infty_0\H$ of the small phase space $\H$.

Introduce functions $\<\<\sigma_{k_1,a_1}\dots\sigma_{k_n,a_n}\>\>^X_g$ on
the large phase space by the formula
\begin{multline*}
\<\<\sigma_{k_1,a_1}\dots\sigma_{k_n,a_n}\>\>^X_g = \sum_{\beta\in
H_2(X,\Z)} q^\beta \sum_{N=0}^\infty \frac{1}{N!}
\sum_{\substack{ b_1\dots b_N \\ \ell_1\dots\ell_N }}
  t^{b_1}_{\ell_1} \dots t^{b_N}_{\ell_N} \\
  \shoveleft{\int_{[\Mbar_{g,n,N}(X,\beta)]^\virt}
  \psi_1^{k_1}\dots\psi_n^{k_n}\Psi_{n+1}^{\ell_1}\dots\Psi_{n+N}^{\ell_N}} \\
\ev_1^*\gamma_{a_1} \dots \ev_n^*\gamma_{a_n}
\ev_{n+1}^*\gamma_{b_1} \dots \ev_{n+N}^*\gamma_{b_N} .
\end{multline*}
\begin{theorem} \label{main}
If $k_i>0$, then $\<\<\sigma_{k_1,a_1}\dots\sigma_{k_n,a_n}\>\>^X_g =
\D_{k_1,a_1}\dots\D_{k_n,a_n}\PHI^X_g$.
\end{theorem}
\begin{proof}
Form the generating function
\begin{multline*}
\begin{aligned}
\PHI^X_{g|a_1\dots a_n}(z_1,\dots,z_n) &= \sum_{k_1,\dots,k_n=0}^\infty
z_1^{k_1} \dots z_n^{k_n} \<\<\sigma_{k_1,a_1}\dots\sigma_{k_n,a_n}\>\>^X_g
\\
&= \sum_{\beta\in H_2(X,\Z)} q^\beta \sum_{N=0}^\infty \frac{1}{N!}
\sum_{\substack{b_1\dots b_N \\ \ell_1\dots\ell_N}} t^{b_1}_{\ell_1} \dots
t^{b_N}_{\ell_N} \end{aligned} \\
\int_{[\Mbar_{g,n,N}(X,\beta)]^\virt}
\frac{\Psi_{n+1}^{\ell_1}\dots\Psi_{n+N}^{\ell_N}
\ev_1^*\gamma_{a_1}\dots\ev_n^*\gamma_{a_n}
\ev_{n+1}^*\gamma_{b_1}\dots\ev_{n+N}^*\gamma_{b_N}}
{(1-z_1\psi_1)\dots(1-z_n\psi_n)} .
\end{multline*}
Corollary \ref{Stab} implies that
$$
\Theta_{a_1}^{b_1}(z_1) \dots \Theta_{a_n}^{b_n}(z_n)
\PHI^X_{g|b_1\dots b_n}(z_1,\dots,z_n) =
\d_{a_1}(z_1)\dots\d_{a_n}(z_n)\PHI^X_g .
$$
By the relation $\D_a(z)=\Theta_a^b(z)\d_b(z)$, it follows that
$$
\Theta_{a_2}^{b_2}(z_2) \dots \Theta_{a_n}^{b_n}(z_n)
\PHI^X_{g|a_1b_2\dots b_n}(z_1,\dots,z_n) =
\D_{a_1}(z_1)\d_{a_2}(z_2)\dots\d_{a_n}(z_n)\PHI^X_g .
$$
Since
$$
\D_{a_1}(z_1) \* \Theta_{a_2}^{b_2}(z_2) = \Theta_{a_2}^{b_2}(z_2) \*
\bigl( \D_{a_1}(z_1) + z_2 \A{a_1} \bigr) ,
$$
we see that
\begin{multline*}
\Theta_{a_3}^{b_3}(z_3) \dots \Theta_{a_n}^{b_n}(z_n)
\PHI^X_{g|a_1a_2b_3\dots b_n}(z_1,\dots,z_n) \\ =
\bigl(\D_{a_1}(z_1)+z_2\A{a_1}\bigr)\D_{a_2}(z_2)
\d_{a_3}(z_3)\dots\d_{a_n}(z_n)\PHI^X_g . 
\end{multline*}
Continuing in the same vein, we see that
\begin{multline*}
\PHI^X_{g|a_1\dots a_n}(z_1,\dots,z_n) \\ = \bigl( \D_{a_1}(z_1) +
(z_2+\dots+z_n) \A{a_1} \bigr) \bigl( \D_{a_{n-1}}(z_{n-1}) + z_n
\A{a_{n-1}} \bigr) \D_{a_n}(z_n) \PHI^X_g .
\end{multline*}
The proposition follows on expanding in the parameters $z_i$.
\end{proof}

This theorem immediately implies Theorem \ref{Fg}. Let $g>0$, so that
the moduli space $\Mbar_{g,n}$ has dimension $3g-3+n$; it follows that
when $k_1+\dots+k_n>3g-3+n$,
$$
\D_{k_1}\dots\D_{k_n}\PHI^X_g=\<\<\D_{k_1,a_1}\dots\D_{k_n,a_n}\>\>^X_g=0 .
$$
Hence, by Proposition \ref{vanish}, we see
that when $k_1+\dots+k_n>3g-3+n$,
$$
\d_{k_1,a_1}\dots\d_{k_n,a_n}\PHI^X_g=0 .
$$

Theorem \ref{main} gives rise to a simple procedure for translating
identities among cycles on Deligne-Mumford moduli spaces into
differential equations among Gromov-Witten potentials. For example,
the identity on $\Mbar_{1,1}$
\begin{equation} \label{trr1}
\psi_1 \cap [\Mbar_{1,1}] = \frac{1}{24}
\begin{picture}(20,25)(70,775)
\put( 80,788){\line( 0,-1){23}}
\put( 80,790){\circle{5}}
\end{picture} \\
\end{equation}
\vskip0.2cm
\noindent where the vertices denoted by open circles are of genus $1$,
translates into the topological recursion relation
$$
\<\<\D_{1,a}\>\>_1^X = \tfrac{1}{24} \<\<\CO_a\CO_b\CO^b\>\>_0^X .
$$
Here, $\CO_a$ is an abbreviation for $\D_{0,a}=\d_{0,a}$. This is the
case $k=1$ of the following equation:
\begin{equation} \label{TRR1}
\D_{k,a}\PHI_1 = \begin{cases} \tfrac{1}{24} \, \D_{0,a} \Tr(\M) & k=1 ,
\\ 0 & k>1 .
\end{cases}
\end{equation}
(The cases $k>1$ reflect the fact that $\Mbar_{1,1}$ has dimension
$1$.)
\begin{lemma}
The equation \eqref{TRR1} has a particular solution $\tfrac{1}{24}
\log\det(\X) \in \CO_1$.
\end{lemma}
\begin{proof}
Clearly, $\D_{k,a}\log\det(\X)=0$ for $k>1$, while
\begin{align*}
\D_{1,a}\log\det(\X) &= \Tr\bigl( \X^{-1} (\X^2){}_a^b \p_{1,b}\X \bigr)
= \Tr\bigl( \X^{-1} (\X^2){}_a^b \A{b} \bigr) \\
&= \Tr\bigl( \X\A{a} \bigr) = \D_{0,a} 
\Tr(\M) .
\qedhere\end{align*}
\end{proof}

Let $\GG = \PHI_1 - \tfrac{1}{24} \log\det(\X)$; we see that
$\D_{k,a}\GG=0$ for all $k>0$; hence, by Corollary \ref{vanish}, we
recover a result of Dijkgraaf and Witten~\cite{DW}: there is a
function $\GG\in\CO_0$ such that
$$
\PHI_1 = \GG + \tfrac{1}{24}\log\det(\X) .
$$
The dilaton equation $\Dil\PHI_1=\tfrac{1}{24}\chi(X)$ follows
automatically from the formulas $\Dil\GG=0$ and
$$
\Dil\log\det(\X)=\Tr(\X^{-1}\Dil\X)=\chi .
$$

To give a second example, the identity on $\Mbar_{2,1}$
\begin{equation*}
\psi_1^2 \cap [\Mbar_{2,1}] = \frac{7}{5}
\begin{picture}(40,35)(60,775)
\put( 80,780){\line( 1, 1){ 13}}
\put( 80,780){\circle*{4}}
\put( 80,780){\line(-1, 1){ 13}}
\put( 80,780){\line( 0,-1){ 18}}
\put( 66,796){\circle{5}}
\put( 94,796){\circle{5}}
\end{picture}
+ \frac{1}{5}
\begin{picture}(30,35)(65,757)
\put( 80,760){\circle*{4}}
\put( 80,782){\circle{5}}
\put( 80,770){\circle{20}}
\put( 80,760){\line(0,-1){18}}
\end{picture}
- \frac{1}{120}
\begin{picture}(30,35)(65,757)
\put( 80,776){\circle{20}}
\put( 80,766){\line( 0,-1){9}}
\put( 80,755){\circle{5}}
\put( 80,753){\line( 0,-1){11}}
\put( 80,766){\circle*{4}}
\end{picture}
+ \frac{13}{120}
\begin{picture}(30,35)(65,757)
\put( 80,770){\circle{20}}
\put( 80,760){\circle*{4}}
\put( 80,760){\line( 3,-4){10}}
\put( 80,760){\line(-3,-4){12}}
\put( 91,745){\circle{5}}
\end{picture}
+ \frac{1}{120}
\begin{picture}(50,35)(55,760)
\put( 70,765){\circle{20}}
\put( 90,765){\circle{20}}
\put( 80,765){\circle*{4}}
\put( 80,765){\line(0,-1){20}}
\end{picture} \\
\end{equation*}
\vskip0.3cm
\noindent translates into the topological recursion relation
\begin{multline*}
\<\<\D_{2,a}\>\>_2^X = \tfrac{7}{5} \<\<\CO_a\CO^b\CO^c\>\>_0^X
\<\<\CO_b\>\>_1^X \<\<\CO_c\>\>_1^X + \tfrac{1}{5}
\<\<\CO_a\CO^b\CO^c\>\>_0^X \<\<\CO_b\CO_c\>\>_1^X \\
- \tfrac{1}{120}
\<\<\CO_a\CO_b\>\>_1^X \<\<\CO^b\CO^c\CO_c\>\>_0^X + \tfrac{13}{120}
\<\<\CO_a\CO_b\CO_c\CO^c\>\>_0^X \<\<\CO^b\>\>_0^X \\ + \tfrac{1}{120}
\<\<\CO_a\CO_b\CO^b\CO_c\CO^c\>\>_0^X .
\end{multline*}

\section{Virasoro constraints on the jet-space of a Frobenius manifold}

\subsection{Frobenius manifolds}
A Frobenius manifold is a conformal Dubrovin manifold
$$
(M,\nabla^z,e,E,r)
$$
together with a non-degenerate symmetric bilinear form $\eta$ on
the tangent bundle (also know as a pseudo-Rimannian metric) such that
\begin{enumerate}
\item $\nabla$ is the Levi-Civit\`a connection associated to $\eta$;
\item $\eta(X\circ Y,Z)=\eta(X,Y\circ Z)$, or equivalently,
$\A{X}^*=\A{X}$;
\item $E\eta(X,Y)=\eta([E,X],Y)+\eta(X,[E,Y])+(2-r)\eta(X,Y)$, that is,
$\CL_E\eta=(2-r)\eta$.
\end{enumerate}

On a Frobenius manifold, the endomorphism $\mu$ of the cotangent bundle is
skew-adjoint:
\begin{multline*}
\eta(\mu^*(X),Y) + \eta(X,\mu^*(Y)) = (2-r) \eta(X,Y) + \eta(\nabla_XE,Y) +
\eta(X,\nabla_YE) \\ = (2-r) \eta(X,Y) + \eta(\nabla_EX,Y) +
\eta(X,\nabla_EY) - \eta([E,X],Y) - \eta(X,[E,Y]) = 0 .
\end{multline*}
Here, we have used that $\nabla$ is torsion-free and compatible with the
metric $\eta$.

If $A$ is a linear operator on the bundle $\CH=T^*M[z,z^{-1}]$, denote by
$A^+$ its adjoint with respect to the symmetric bilinear form
$$
\Eta( f , g ) = \Res_{z=0} \eta(f(-z),g(z)) \frac{dz}{z} .
$$
For example, $\Theta^+(z)=\Theta^*(-z)$, $\p_z^+=\p+z^{-1}$ and $z^+=-z$.
Since $\mu^*=-\mu$ and $\U^+=\U$, it follows that
$$
\delta_z^+ = ( \p_z + z^{-1} ) - z^{-1} ( \mu^* + \half ) + \U^*
= ( \p_z + z^{-1} ) + z^{-1} ( \mu - \half ) + \U = \delta_z .
$$

If $\Theta(z)$ is a fundamental solution, recall that the endomorphism
$R(z)$ of $\CH=T^*M[z]$ is defined by the formula
$$
\p_z\Theta(z) + z^{-1} R(z)\Theta(z) = z^{-1} \bigl( \mu + \half \bigr)
\Theta(z) + \U \Theta(z) .
$$
\begin{definition}
A fundamental solution $\Theta(z)$ on a Frobenius manifold $M$ is symmetric
if it is conformal and $\Theta^+(z)\Theta(z)=\Id$.
\end{definition}

It follows from the condition $\Theta^+(z)\Theta(z)=\Id$ that
$R(z)+R^+(z)=0$, in other words, $R_k^*=(-1)^{k+1}R_k$ for all $k>0$.

\begin{proposition}
There exists a symmetric fundamental solution on any Frobenius manifold,
\end{proposition}
\begin{proof}
Let $\Theta(z)$ be a conformal fundamental solution. Since
$\nabla\*\Theta(z)=\Theta(z)\*\nabla^z$, we have
\begin{align*}
0 &= \eta(\alpha(-z),\Theta(z)\*\nabla^z\beta(z)) -
\eta(\alpha(-z),\nabla\*\Theta(z)\beta(z)) \\
&= \eta(\Theta^*(z)\alpha(-z),\nabla^z\beta(z)) -
d\eta(\alpha(-z),\Theta(z)\beta(z))
+ \eta(\nabla\alpha(-z),\Theta(z)\beta(z)) \\
&= - \eta(\nabla^{-z}\*\Theta^*(z)\alpha(-z),\beta(z))
+ \eta(\Theta(z)\*\nabla\alpha(-z),\beta(z)) .
\end{align*}
This shows that $\nabla^{-z}\*\Theta^*(z)=\Theta^*(z)\*\nabla$; conjugating
by $\Theta^*(z)$, we obtain the equation
$$
\nabla\*\Theta^+(-z) = \Theta^+(z)\*\nabla^{-z} .
$$
Replacing $z$ by $-z$, we see that $\Theta^+(z)^{-1}$ is also a fundamental
solution.

It follows that $\rho(z)=\Theta^+(z)\Theta(z)$ is a flat section of
$\Gamma(M,\End(TM))\[z\]$. Replacing $\Theta(z)$ by
$\tTheta(z)=\Theta(z)\rho(z)^{-1/2}$, we obtain a new fundamental
solution. Since $\rho^+(z)=\rho(z)$, we see that
$$
\tTheta^+(z)\tTheta(z) = \rho^+(z)^{-1/2}\Theta^+(z)\Theta(z)\rho(z)^{-1/2}
= \Id ,
$$
hence $\tTheta(z)$ is symmetric.
\end{proof}

In the remainder of this section, we only consider symmetric
fundamental solutions.

\subsection{Gromov-Witten invariants and Frobenius manifolds}

In Section 4, we showed that the small phase space $\H$ associated to
the genus $0$ Gromov-Witten invariants of a projective manifold $X$ is
conformal Dubrovin manifold. It also carries a flat metric $\eta$;
with respect to this metric, it is a Frobenius manifold. The equation
$\CL_E\eta=(2-r)\eta$ is equivalent to
$(p_a+p_b)\eta_{ab}=r\eta_{ab}$, which is a basic property of the
Poincar\'e form of a projective manifold of dimension $r$.
\begin{proposition}
The conformal fundamental solution
$$
\Theta_a^b(z) = \delta_a^b + \eta^{bc} \sum_{k=0}^\infty z^{k+1}
\<\<\tau_{k,a}\tau_{0,c}\>\>^X_0
$$
on the small phase space is symmetric.
\end{proposition}
\begin{proof}
Let $\rho(z)=\Theta^+(z)\Theta(z)$. We have
\begin{align*}
\p_a\rho^+(z) &= \p_a\Theta^+(z)\Theta(z) + \Theta^+(z)\p_a\Theta(z) \\
&= -z\Theta^+(z) \A{a} \Theta(z) + z \Theta^+(z) \A{a} \Theta(z) = 0 .
\end{align*}
Thus $\rho(z)$ is constant on $\H$. It remains to show that it equals $\Id$
at the origin $0\in\H$. We may write
$$
\rho(z)(0) = \Id + \sum_{k=0}^\infty z^{k+1}
\sum_{\beta\ne0} q^\beta \rho_{k,\beta} .
$$
For $\om\in H^2(X,\C)$, the divisor equation implies that
$$
\tint_\beta \om \, \rho_{k,\beta} = [R(\om),\rho_{k-1,\beta}] ,
$$
from which it follows that $\rho_{k,\beta}=0$ for $\beta\ne0$. On the other
hand, in the limit $q\to0$ in the Novikov ring, the value of $\Theta_k$ at
the origin of $\H$ is an integral over $\Mbar_{0,2}(X,0)$, hence vanishes,
since $\Mbar_{0,2}(X,0)$ is empty.
\end{proof}

\subsection{The genus $0$ potential on the large phase space of a Frobenius
manifold}
On a Frobenius manifold, we may define a power series $\PHI_0$ on the large
phase space by integrating the power series $\g(z)$:
$$
\d_{k,a}\PHI_0 = (-1)^k \eta_{ab} \g_{-k-1}^b .
$$
Denote $\eta_{bc}\Om_{j,k|a}^c$ by $\Om_{j,k|a,b}$; taking another
derivative, we have
\begin{equation} \label{ddF}
\d_{j,a}\d_{k,b}\PHI_0 = (-1)^k \Om_{j,k|a,b} ,
\end{equation}
This determines $\PHI_0$ up to an affine function on the large phase space,
which may be fixed by the dilaton equation
$$
\Dil\PHI_0=-2\PHI_0 .
$$
The following formula $\PHI_0$ is due to Dubrovin \cite{D}.
\begin{proposition}
$
\PHI_0 = \half \sum_{\ell,m=0}^\infty (-1)^m \, \g_\ell^p \, \g_m^q \,
\Om_{\ell,m|p,q}
$
\end{proposition}
\begin{proof}
Since the fundamental solution $\Theta(z)$ is symmetric, we see that
$$
(y-z) \Om^*(y,z) = \Theta^*(z)^{-1}\Theta^*(y) - \Id
= \Theta(-z)\Theta^{-1}(-y) - \Id = (y-z) \Om(-z,-y) ,
$$
in other words, $\Om_{j,k|a,b} = (-1)^{j+k} \Om_{k,j|b,a}$. It follows by
Lemma \ref{integrable} that
$$
\d_{k,b} \Om_{\ell,m|p,q} = (-1)^{\ell+m} \d_{k,b} \Om_{m,\ell|q,p}
= (-1)^{\ell+m} \d_{m,q} \Om_{k,\ell|b,p}
= (-1)^{k+m} \d_{m,q} \Om_{\ell,k|p,b} .
$$

We have
\begin{align*}
\d_{k,b} \PHI_0 &= \half \sum_{\ell,m} (-1)^m \g_\ell^p \, \g_m^q \,
\d_{k,b} \Om_{\ell,m|p,q} \\ &\quad + \half (-1)^k \sum_\ell \g_\ell^p \,
\Om_{\ell,k|p,b} + \half \sum_m (-1)^m \g_m^q \, \Om_{k,m|b,q} \\
&= \half (-1)^k \sum_{\ell,m} \g_\ell^p \, \g_m^q \, \d_{m,q}
\Om_{\ell,k|p,b} \\ &\quad + \half (-1)^k \sum_\ell \g_\ell^p \,
\Om_{\ell,k|p,b} + \half \sum_m (-1)^m \g_m^q \, \Om_{k,m|b,q} \\
&= \half (-1)^k \sum_\ell \g_\ell^p \, ( 1 - \Dil ) \Om_{\ell,k|p,b}
+ \half \sum_m (-1)^m \g_m^q \, \Om_{k,m|b,q} \\
&= \half (-1)^k \sum_\ell \g_\ell^p \, \Om_{\ell,k|p,b} + \half \sum_m
(-1)^m \g_m^q \, \Om_{k,m|b,q} ,
\end{align*}
since $\Dil\Om(y,z)=0$. It follows that
\begin{align*}
\d_{j,a} \d_{k,b} \PHI_0 &= \half (-1)^k \Om_{j,k|a,b}
+ \half (-1)^k \sum_\ell \g_\ell^p \, \d_{j,a} \Om_{\ell,k|p,b} \\ &\quad
+ \half (-1)^j \Om_{k,j|b,a} + \half \sum_m (-1)^m \g_m^q \, \d_{j,a}
\Om_{k,m|b,q} \\
&= \half (-1)^k \Om_{j,k|a,b} + \half (-1)^k \sum_\ell \g_\ell^p \,
\d_{\ell,p} \Om_{j,k|a,b} \\ &\quad+ \half (-1)^j \Om_{k,j|b,a} +
\half (-1)^j \sum_m \g_m^q \, \d_{m,q} \Om_{k,j|b,a} \\
&= \half (-1)^k (1-\Dil)\Om_{j,k|a,b} + \half (-1)^j (1-\Dil) \Om_{k,j|b,a}
\\
&= \half (-1)^k \Om_{j,k|a,b} + \half (-1)^j \Om_{k,j|b,a} = (-1)^k
\Om_{j,k|a,b} .
\end{align*}
This shows that $\PHI_0$ satisfies \eqref{ddF}.
\end{proof}

\subsection{The Virasoro operators on the large phase space of a Frobenius
manifold}
Let $\Phi(z)$ be the generating function (or \emph{free field}) whose
coefficients are linear differential operators acting on functions on the
large phase space
$$
\Phi_a(z) = \hbar \, \sum_{k=0}^\infty z^k \d_{k,a}
+ \sum_{k=0}^\infty (-z)^{-k-1} \eta_{ab} \g_k^b .
$$
The conjugate of $\Phi(z)$ by $Z_0=e^{\PHI_0/\hbar}$ is given by the
formula
$$
Z_0^{-1} \* \Phi_a(z) \* Z_0 = \hbar \, \d_a(z) + \eta_{ab} \g^b(-z) .
$$
The normal ordering of quadratic expressions in the field $\Phi(z)$ is defined
by the formulas
\begin{align*}
\normal{\g^a_k\,\g^b_\ell} &= \g^a_k\,\g^b_\ell , &
\normal{\g^a_k\,\d_{\ell,b}} &= \g^a_k\,\d_{\ell,b} , &
\normal{\d_{k,b}\,\g^a_\ell} &= \g^a_\ell\,\d_{k,b} , &
\normal{\d_{k,a}\,\d_{\ell,b}} &= \d_{k,a}\,\d_{\ell,b} .
\end{align*}
If $k\ge-1$, let $L_k$ be the second-order differential operator on the
large phase space
$$
L_k = \tfrac{1}{2\hbar} \Res_{z=0} \normal{\eta( \Phi(-z) ,
z\tdelta_z^{k+1}\Phi(z) )} - \tfrac{1}{4} \delta_{k,0} \Tr\bigl( \mu^2 -
\tfrac{1}{4} \bigr) .
$$
These operators were introduced in the theory of Gromov-Witten invariants
by Eguchi, Hori and Xiong \cite{EHX}, who proved the following result.
\begin{proposition}
The operators $L_k$ satisfy the commutation relations
$$[L_j,L_k]=(j-k) L_{j+k}.$$
\end{proposition}
\begin{proof}
We give an outline of the proof; see \cite{virasoro}, Section 2, for the
complete details.

Using the canonical commutation relations among the coefficients of
$\Phi(z)$, the formula $[L_k,\Phi(z)]+z\tdelta_z^{k+1}\Phi(z)=0$ is easily
proved. It follows that
\begin{align*}
[[L_j,L_k],\Phi(z)] &= [L_j,[L_k,\Phi(z)]] - [L_k,[L_j,\Phi(z)]] \\
&= - [L_j,z\tdelta_z^{k+1}\Phi(z)] + [L_k,z\tdelta_z^{j+1}\Phi(z)] \\
&= - z\tdelta_z^{k+1} [L_j,\Phi(z)] + z\tdelta_z^{j+1} [L_k,\Phi(z)] \\
&= [z\tdelta_z^{k+1},z\tdelta_z^{j+1}]\Phi(z) = (k-j)
\tdelta_z^{j+k+1}\Phi(z) ,
\end{align*}
and hence that $[L_j,L_k]-(j-k)L_{j+k}$ commutes with $\Phi(z)$.

Any differential operator commuting with $\Phi(z)$ must lie in the centre of
the algebra of differential operators, and hence
$$
[L_j,L_k] = (j-k)L_{j+k} + c(j,k)1
$$
for some two-cocyle $c(j,k)$. Evaluating both sides of this
equation at the basepoint of the large phase space, the result follows.
\end{proof}

In the following theorem, we show how the Virasoro operators $L_k$
give rise to intrinsic differential operators on the jet-space of the
Frobenius manifold. The explicit formula for $\h_k$ is equivalent to a
formula of Liu (\cite{Liu}, Theorem 4.4).
\begin{theorem} \label{transfer}
We have
$$
L_k\*Z_0 = Z_0 \* \bigl( \hbar\,\Delta_k + \L_k + \h_k \bigr) ,
$$
where $\Delta_k = \half \Res_{z=0}\eta\bigl( \D(-z) , z\delta_z^{k+1} \D(z)
\bigr)$, and
$$
\h_k = - \tfrac{1}{4} \sum_{\ell=0}^k \Tr\bigl( \bigl(\mu-\half\bigr)\U^\ell
\bigl(\mu+\half\bigr)\U^{k-\ell} \bigr) .
$$
\end{theorem}
\begin{proof}
We have
\begin{multline*}
Z_0^{-1} \* \tfrac{1}{2\hbar} \normal{\eta( \Phi(-z) ,
z\tdelta_z^{k+1}\Phi(z) )} \* Z_0
= \half \hbar \, \eta\bigl( \d(-z) , z\tdelta_z^{k+1}\d(z) \bigr) \\
\begin{aligned}
&+ \half \bigl\< \g(z) , z\tdelta_z^{k+1}\d(z) \bigr\>
+ \half \bigl\< z\tdelta_z^{k+1}\g_-(-z) , \d(-z) \bigr\>
+ \half \bigl\< \d(-z) , z\tdelta_z^{k+1}\g_+(-z) \bigr\> \\
&+ \tfrac{1}{2\hbar} \eta\bigl(\g(z),z\tdelta_z^{k+1}\g(-z)\bigr) .
\end{aligned}
\end{multline*}
On taking the residue, the first term on the right-hand side yields
\begin{multline*}
\half \Res_{z=0} \eta\bigl( \d(-z) , z\tdelta_z^{k+1} \d(z) \bigr) \\ =
\half \Res_{z=0} \eta\bigl( \Theta^{-1}(-z)\D(-z) , \Theta^{-1}(z)
z\delta_z^{k+1} \D(z) \bigr) = \Delta_k .
\end{multline*}
The residue of the last term on the right-hand side vanishes by Theorem
\ref{G} (this is the genus $0$ Virasoro constraint), since
\begin{align*}
\eta\bigl(\g(z),z\tdelta_z^{k+1}\g(-z)\bigr)
&= \eta\bigl(\Theta(z)\g(z),\Theta^{-1}(z)z\tdelta_z^{k+1}\g(-z)\bigr)
\\ &= \eta\bigl( \G(z),z\delta_z^{k+1}\G(-z) \bigr) = O(z^{-3}) .
\end{align*}
It follows that
\begin{multline*}
Z_0^{-1} \* L_k \* Z_0 - \hbar \, \Delta_k + \tfrac{1}{4} \delta_{k,0}
\Tr\bigl( \mu^2 - \tfrac{1}{4} \bigr) \\
\begin{aligned}
&= \half \Res_{z=0} \bigl\< \g(z) , z\tdelta_z^{k+1}\d(z) \bigr\> \\
  &\quad+ \half \Res_{z=0} \Bigl( \bigl\< z\tdelta_z^{k+1}\g_-(-z) ,
  \d(-z) \bigr\> + \bigl\< \d(-z) , z\tdelta_z^{k+1}\g_+(-z) \bigr\>
  \Bigr) \\
&= \half \Res_{z=0} \Bigl( \bigl\< \g(z) , z\tdelta_z^{k+1}\d(z) \bigr\> +
\bigl\< \g_-(z) , z\tdelta_z^{k+1}\d(z) \bigr\> + \bigl\<
z\tdelta_z^{k+1}\d(z) , \g_+(z) \bigr\> \Bigr) \\
&= \Res_{z=0} \bigl\< \g(z) , z\tdelta_z^{k+1}\d(z) \bigr\> -
\half \Res_{z=0} \Tr\bigl( y\tdelta_y^{k+1}\Om(y,z) \bigr){}_{y=z} ,
\end{aligned}
\end{multline*}
since $\d(y)\g_+(z)=-\Om(y,z)$. It remains to show that
\begin{multline*}
\Res_{z=0} \Tr\bigl( y\tdelta_y^{k+1}\Om(y,z) \bigr){}_{y=z} \\
= \begin{cases}
0 , & k\le0, \\
\displaystyle - \half \sum_{\ell=0}^k \Tr\bigl( \bigl(\mu-\half\bigr)\U^\ell
\bigl(\mu+\half\bigr)\U^{k-\ell} \bigr) , & k>0 .
\end{cases}
\end{multline*}
Applying the operator $y\tdelta_y^{k+2}$ to \eqref{Om}, we see that
$$
(k+2) y\tdelta_y^{k+1}\Om(y,z) + (y-z) y\tdelta_y^{k+2} \Om(y,z) =
y\tdelta_y^{k+2} \bigl( \Theta(y)\Theta^{-1}(z) - \Id \bigr) .
$$
Taking the trace, we see that
\begin{multline*}
(k+2) \Tr\bigl( y\tdelta_y^{k+1}\Om(y,z) \bigr) + (y-z) \Tr\bigl(
y\tdelta_y^{k+2} \Om(y,z) \bigr) \\ =
\Tr\bigl( \Theta^{-1}(z) y\tdelta_y^{k+2}\Theta(y) \bigr)
- \Tr\bigl( y\tdelta_y^{k+2}\Id \bigr) .
\end{multline*}
Next, we restrict to the diagonal $y=z$:
\begin{align*}
(k+2) \Tr\bigl( y\tdelta_y^{k+1}\Om(y,z) \bigr)_{y=z} &=
\Tr\bigl( \Theta^{-1}(z) z\tdelta_z^{k+2}\Theta(z) \bigr)
- \Tr\bigl( z\tdelta_z^{k+2}\Id \bigr) \\
&= \Tr\bigl( z\delta_z^{k+2}\Id \bigr) -
\Tr\bigl( z\tdelta_z^{k+2}\Id \bigr) .
\end{align*}
Taking the residue, we see that
\begin{multline*}
(k+2) \Res_{z=0} \Tr\bigl( y\tdelta_y^{k+1}\Om(y,z) \bigr)_{y=z} \\
= \sum_{0\le i+j\le k} \Tr\bigl( \U^i \bigl( \mu-\half \bigr) \U^j \bigl(
\mu + \half \bigr) \U^{k-i-j} \bigr) - \Res_{z=0}\Tr\bigl(
z\tdelta_z^{k+2}\Id \bigr) .
\end{multline*}
Since $[\mu_s,R_k]=kR_k$, we see that
$$
\Tr\bigl( z\tdelta_z^{k+1}\Id \bigr) = z^{-k-1} \Tr\bigl( \bigl( \mu -
\half - k \bigr) \dots \bigl( \mu - \half \bigr) \bigl( \mu + \half \bigr)
\bigr) ,
$$
and hence $\Res_{z=0} \Tr\bigl( z\tdelta_z^{k+1}\Id \bigr) = \delta_{k,0}
\Tr\bigl(\mu^2 - \quart \bigr)$. We have
\begin{multline*}
\sum_{0\le i+j\le k} \Tr\bigl( \U^i \bigl( \mu-\half \bigr) \U^j \bigl(
\mu + \half \bigr) \U^{k-i-j} \bigr) \\
\begin{aligned}
&= \sum_{0\le i+j\le k} \Tr\bigl( \U^i\mu\U^j\mu\U^{k-i-j} \bigr)
- \quart \sum_{0\le i+j\le k} \Tr\bigl( \U^k \bigr) \\
&= \bigl( \tfrac{k}{2} + 1 \bigr) \sum_{\ell=0}^k \Tr\bigl(
\mu\U^\ell\mu\U^{k-\ell} \bigr) - \quart \tbinom{k+2}{2} \Tr\bigl( \U^k
\bigr) \\
&= \bigl( \tfrac{k}{2} + 1 \bigr) \sum_{\ell=0}^k \Tr\bigl( \bigl(
\mu-\half \bigr) \U^\ell \bigl( \mu + \half \bigr) \U^{k-\ell} \bigr) ,
\end{aligned}
\end{multline*}
and the formula for $\h_k$ follows.
\end{proof}

As a corollary of this theorem and the Virasoro relations for $L_k$, we
obtain the following formula of Liu (\cite{Liu}, Section 6).
\begin{corollary}
$\L_j\h_k-\L_k\h_j=(j-k)\h_{j+k}$
\end{corollary}

\subsection{The Virasoro constraints on a Frobenius manifold}
The Virasoro constraints are differential equations among a sequence of
functions $\PHI_g$, $g>0$, on the large phase space of a Frobenius manifold
$M$ such that the dilaton equation holds:
\begin{equation} \label{Dilaton}
\Dil\PHI_g = (2g-2)\PHI_g + \tfrac{1}{24} \delta_{g,1} \chi ,
\end{equation}
were $\chi=\Tr(\Id)$. As originally introduced by Eguchi, Hori and Xiong
\cite{EHX}, in the case where $M$ is the Frobenius manifold associated to
the genus $0$ Gromov-Witten theory of a projective manifold $X$ and
$\PHI_g$ is the genus $g$ Gromov-Witten potentials $X$, these differential
equations are given by the vanishing of functions $z_{k,g}$ on the large
phase space defined as follows:
$$
\sum_{g=0}^\infty \hbar^{g-1} z_{k,g} = \textstyle \exp\bigl( -
\sum_{g=0}^\infty \hbar^{g-1} \PHI_g \bigr) \* L_k \* \exp\bigl(
\sum_{g=0}^\infty \hbar^{g-1} \PHI_g \bigr) .
$$
Using Theorem \ref{transfer}, we may write the vanishing of $z_{k,g}$
intrinsically on the jet-space of the Frobenius manifold $M$: in fact, we
have
\begin{multline*}
\sum_{g=1}^\infty \hbar^{g-1} z_{k,g} = \textstyle \exp\bigl( -
\sum_{g=1}^\infty \hbar^{g-1} \PHI_g \bigr) \* \hbar\,\Delta_k \*
\exp\bigl( \sum_{g=1}^\infty \hbar^{g-1} \PHI_g \bigr) \\
+ \sum_{g=1}^\infty \hbar^{g-1} \L_k\PHI_g + \h_k .
\end{multline*}
This proves the following theorem.
\begin{theorem} \label{Virasoro}
For genus $g>1$, the Virasoro constraint $z_{k,g}=0$ is equivalent to
the equation
$$
\L_k\PHI_g + \Delta_k\PHI_{g-1} + \half \sum_{i=1}^{g-1} \Res_{z=0}
\eta\bigl( \D(-z)\PHI_i , z\delta_z^{k+1}\D(z)\PHI_{g-i} \bigr) = 0 ,
$$
while the Virasoro constraint $z_{k,1}=0$ is equivalent to
the equation
$$
\L_k\PHI_1 = \tfrac{1}{4} \sum_{\ell=0}^k \Tr\bigl(
\bigl(\mu-\half\bigr)\U^\ell\bigl(\mu+\half\bigr)\U^{k-\ell} \bigr) .
$$
\end{theorem}

Using Proposition \ref{EkX}, we see that the genus $1$ Virasoro constraint
$z_{k,1}=0$ is equivalent to the equation
\begin{equation} \label{V1}
\L_k\GG = \sum_{\ell=0}^k \Bigl( \tfrac{1}{4}
\Tr\bigl(\mu\U^\ell\mu\U^{k-i} \bigr) + \tfrac{1}{24} \bigl< e ,
\U^\ell\mu\U^{k-\ell} \om \bigr\> \Bigr)
\end{equation}
for $\GG=\PHI_1^X-\frac{1}{24}\log\det(\X)$. This equation was first
proved in the special case where $M$ is a semisimple Frobenius
manifold by Dubrovin and Zhang \cite{DZ}, and extended to general
Frobenius manifolds by Liu \cite{Liu}.

Let $M$ be a Frobenius manifold. There is a second-order differential operator
$$
\Psi : \Gamma(M,\CO) \to \Gamma(M,S^4(T^*M)) ,
$$
given by the explicit formula
\begin{multline} \label{elliptic}
\Psi(X_1,X_2,X_3,X_4) = \tfrac{1}{24} \sum_{\pi\in S_4}
\sum_{a_1,a_2,a_3,a_4}
X_{\pi_1}^{a_1}X_{\pi_2}^{a_2}X_{\pi_3}^{a_3}X_{\pi_4}^{a_4} \\
\bigl( 3\,\AA^\mu_{a_1a_2}\,\AA^\nu_{a_3a_4}\,\p_\mu\p_\nu -
4\,\AA^\mu_{a_1a_2}\,\AA^\nu_{a_3\mu} \, \p_{a_4}\p_\nu + 2 \,
\p_{a_4}\AA^\mu_{a_2a_3}\,\AA^\nu_{a_1\mu}\p_\nu - \AA^\mu_{a_1a_2} \,
\p_\mu\AA^\nu_{a_3a_4}\p_\nu \\
+ \tfrac16 \p_{a_3}\AA^\mu_{a_1a_2}\,\p_{a_4}\AA^\nu_{\mu\nu}
+ \tfrac1{24} \p_{a_3}\p_{a_4}\AA^\mu_{a_1a_2}\,\AA^\nu_{\mu\nu}
- \tfrac14 \p_\nu\AA^\mu_{a_1a_2} \, \p_\mu\AA^\nu_{a_3a_4} \bigr) ,
\end{multline}
such that the generating function $\GG$ of genus $1$ Gromov-Witten
invariants of a compact symplectic manifold $X$ satisfies the equation
$\Psi\GG=0$ (see \cite{elliptic}). Let
$$
y_k = \L_k\GG - \sum_{\ell=0}^k \Bigl( \tfrac{1}{4}
\Tr\bigl(\mu\U^\ell\mu\U^{k-i} \bigr) + \tfrac{1}{24} \bigl< e ,
\U^\ell\mu\U^{k-\ell} \om \bigr\> \Bigr) .
$$
By a lengthy calculation, Liu \cite{Liu} proves the formula
$$
y_k = \half \sum_{i=1}^k \Psi(E^{\circ k-i+1},E^{\circ i},E,E)\GG -
\tfrac{k+1}{2} \L_{k-1}y_1 .
$$
In this way, he obtains the following theorem.
\begin{theorem}
  If $\Psi\GG=0$ and $\GG$ satisfies \eqref{V1} for $k=1$, then $\GG$
  satisfies \eqref{V1} for all $k>1$.
\end{theorem}

\subsection*{Relationship to prior results}
The results of Sections 5.1 and 5.2 are due to Dubrovin. The formula for
the Virasoro operators on the large phase space in terms of free fields
(Section~5.3) is a modification of the formalism of Section 2 of
\cite{virasoro}. (Related formulas have been used by Givental
\cite{givental}.)

\section{The jet-space of a semisimple Dubrovin manifold}

\subsection{Semisimple Dubrovin manifolds}
A Dubrovin manifold is \textbf{semisimple} if the subset $M_0\subset
M$ on which the commutative algebra $(T_xM,\circ)$ is semisimple is
dense in $M$. These manifolds were extensively studied by Dubrovin
\cite{D}, who proved the following theorem.
\begin{theorem}
  Around each point semisimple point in a semisimple Dubrovin
  manifold, there are coordinates $v^i$ such that the vector fields
  $\pi_i=\p/\p v^i$ satisfy
$$
\pi_i \circ \pi_j = \delta_{ij} \pi_i .
$$
\end{theorem}
\begin{proof}
  Locally in $M_0$, there is a frame $\{\pi_1,\dots,\pi_n\}$ of the
  tangent bundle $TM$ such that $\pi_i \circ \pi_j = \delta_{ij}
  \pi_i$. We must show that the vector fields $\pi_i$ satisfy
  $[\pi_i,\pi_j]=0$. If $i$ and $j$ are distinct, the equation
\begin{equation} \label{ijk}
\nabla_i(\pi_j\circ\pi_k) - \pi_j\circ(\nabla_i\pi_k) -
\nabla_j(\pi_i\circ\pi_k) + \pi_i\circ(\nabla_j\pi_k) =
[\pi_i,\pi_j]\circ\pi_k
\end{equation}
is a consequence of condition 4) of Proposition \ref{dubrovin}. If $k$
is not equal to $i$ or $j$, \eqref{ijk} implies that
$[\pi_i,\pi_j]\circ\pi_k$ vanishes. On the other hand, if $k$ equals
$i$, we see that
$$
- \pi_j\circ\nabla_i\pi_i - \nabla_j\pi_i + \pi_i\circ(\nabla_j\pi_i) =
[\pi_i,\pi_j]\circ\pi_i ,
$$
from which we conclude that $[\pi_i,\pi_j]\circ\pi_i$ vanishes.
This shows that $[\pi_i,\pi_j]=0$. It follows that there are
coordinates $v^i$ defined locally on $M_0$ such that $\pi_i=\p/\p
v^i$.
\end{proof}

The canonical coordinates $(v^i)$ are seen to be defined up to
addition of a constant and permutation. Let $J=(J^a_i)$ be the
Jacobian
$$
J_i^a = \p u^a/\p v^i ,
$$
and let $J^{-1}=(J^i_a)$ be its inverse
$$
J_a^i = \p v^i/\p u^a .
$$
Denote by $v^i_n$ the functions
$$
v^i_n = \p^nv^i = \p^{n-1} \bigl( J^i_a u^a_1 \bigr)
= \sum_{k=1}^n \tbinom{n-1}{k-1} \bigl(\p^{n-k}J^i_a\bigr) u^a_k .
$$
The identity vector field on $M_0$ is given in canonical
coordinates by the formula
$$
e = \sum_i \pi_i .
$$

\subsection{Conformal semisimple Dubrovin manifolds}
Let $(M,\nabla^z,e,E,r)$ be a semisimple Dubrovin manifold which is
conformal. Taking $X=Y=\pi_i$ in \eqref{E}, we see that
$$
[E,\pi_i] = 2\pi_i \circ [E,\pi_i] + \pi_i ,
$$
hence $[E,\pi_i]=-\pi_i$. This shows that, after shifting the
canonical coordinates $v_i$ by constants $c^i$, the Euler vector field
on $M_0$ is given by the formula
$$
E = \sum_i v^i \pi_i .
$$
In other words, the canonical coordinates are the eigenvalues of
$\U$.
\begin{proposition}
  The semisimple locus $M_0\subset M$ contains the \textbf{tame locus}
  $M_{00}$ of points at which the eigenvalues of $\U$ are distinct and
  nonzero.
\end{proposition}
\begin{proof}
  Locally in $M_{00}$, the endormorphism $\U^*$ has eigenvectors
  $\pi_i$ with eigenvalues $v_i$. We have
$$
\U^*(\pi_i\circ\pi_j) = (\U^*\pi_i)\circ\pi_j = v^i\pi_i\circ\pi_j .
$$
Antisymmetrizing in $i$ and $j$ gives $(v^i-v^j)\pi_i\circ\pi_j =
0$, showing that $\pi_i\circ\pi_j=0$ if $i$ and $j$ are distinct.
Since $E\circ\pi_i=v^i\pi_i$ is nonzero, we see that
$\pi_i\circ\pi_i=f_i\pi_i$, where $f_i$ is nowhere vanishing, hence
$(T_xM,\circ)$ is semisimple in $M_{00}$.
\end{proof}

We will restrict attention to conformal semisimple Dubrovin manifolds
satisfying the following condition.
\begin{condition} \label{gamma}
  There is a section $\gamma\in\Gamma(M_0,\End(T^*M))$ such that
$$
\mu=[\gamma,\U] .
$$
\end{condition}

Let $\mu_i^j$ and $\gamma_i^j$ be the matrix elements of the tensors
$\mu$ and $\gamma$ in the frame $(\pi_i)$; Condition \ref{gamma}
amounts to the relation $\mu_i^j = (v^i-v^j)\gamma_i^j$.

\subsection{Semisimple Frobenius manifolds}
A Frobenius manifold is semisimple if the underlying Dubrovin manifold
is semisimple. We now show, borrowing from the exposition of Manin
\cite{M}, that semisimple Frobenius manifolds satisfy Condition
\ref{gamma}.

In Gromov-Witten theory, it is known that the Frobenius manifolds
associated to Grassmannians \cite{Abrams}, to Del Pezzo varieties
\cite{BayerManin}, and to arbitrary products of these, are semisimple.
Another class of examples are the quotients $V/W$ of a Euclidean
vector space $V$ by a Coxeter group $W$ of reflections, considered by
Saito \cite{Saito} (see also Dubrovin \cite{coxeter}).

Let $M$ be a semisimple Frobenius manifold. If $\eps\in\Om^1(M)$ is
the one-form on $M$ defined by the formula $\eps(X)=\eta(e,X)$, then
$$
\eta(X,Y) = \eta(e\circ X,Y) = \eps(X\circ Y) ;
$$
thus $\eta(X,Y)=0$ if $X\circ Y=0$. It follows that the idempotent
eigenvectors $\pi_i$ on $M_0$ are orthogonal. Let
$\eta_i=\eta(\pi_i,\pi_i)$, let $\eta_{ij}=\pi_j(\eta_i)$, and let
$\eta_{ijk}=\pi_k(\eta_{ij})$. We now list the main properties of
these functions; these are known as the Darboux-Egoroff equations.
\begin{lemma} \label{etaij} \mbox{}
\begin{enumerate}
\item $\eta_{ij} = \eta_{ji}$
\item $\eta_{ijk} = \half \bigl( \frac{\eta_{ki}\eta_{ij}}{\eta_i} +
\frac{\eta_{ij}\eta_{jk}}{\eta_j} + \frac{\eta_{jk}\eta_{ki}}{\eta_k}
\bigr)$ when $i$, $j$ and $k$ are distinct
\item $e(\eta_i)=0$
\item $E(\eta_i)=-r\eta_i$
\end{enumerate}
\end{lemma}
\begin{proof}
The coefficients $\Gamma_{ijk}=\eta(\nabla_i\pi_j,\pi_k)$ of the
Levi-Civit\`a connection $\nabla$ are given by the formula
$$
\Gamma_{ijk} = \half \bigl( \delta_{ik} \eta_{ij} +
\delta_{jk} \eta_{ji} - \delta_{ij} \eta_{ik} \bigr) .
$$
Taking the inner product of \eqref{ijk} with $\pi_\ell$, we obtain the
equation
$$
\half ( \delta_{i\ell}\delta_{jk} + \delta_{j\ell}\delta_{ik} ) \bigl(
\eta_{ij} - \eta_{ji} \bigr) = 0 ,
$$
which shows that $\eta_{ij}=\eta_{ji}$.

A lengthy calculation shows that
\begin{align*}
R_{ijk\ell} &= \quart \bigl( \delta_{i\ell} - \delta_{j\ell} \bigr) \Bigl(
\frac{\eta_{ki}\eta_{ij}}{\eta_i} + \frac{\eta_{ij}\eta_{jk}}{\eta_j} +
\frac{\eta_{jk}\eta_{ki}}{\eta_k} - 2 \eta_{ijk} \Bigr) \\
& + \quart \bigl( \delta_{jk} - \delta_{ik} \bigr)
\Bigr( \frac{\eta_{\ell i}\eta_{ij}}{\eta_i} + 
\frac{\eta_{ij}\eta_{j\ell}}{\eta_j} +
\frac{\eta_{j\ell}\eta_{\ell i}}{\eta_\ell} - 2 \eta_{ij\ell} \Bigr) \\
& + \quart \bigl( \delta_{j\ell} \delta_{ik} - \delta_{i\ell} \delta_{jk}
\bigr) \sum_m \frac{\eta_{im}\eta_{mj}}{\eta_m} .
\end{align*}
Part 2) follows, by the flatness of the connection $\nabla$.

Since the vector field $e$ is flat, we have
\begin{align*}
0 &= \eta(\pi_i,\nabla_ie) = \sum_j \eta(\pi_i,\nabla_i\pi_j) = \sum_j
\Gamma_{iji} \\
&= \half \sum_j \bigl( \eta_{ij} + \delta_{ij} \eta_{ji} - \delta_{ij}
\eta_{ii} \bigr) = \half e(\eta_i) .
\end{align*}

Since $\CL_E\eta=(2-r)\eta$, we see that
$$
E(\eta_i) = (\CL_E\eta)(\pi_i,\pi_i) - 2 \eta([E,\pi_i],\pi_i) = (2-r)
\eta(\pi_i,\pi_i) - 2 \eta(\pi_i,\pi_i) ,
$$
which shows that $E(\eta_i)=-r\eta_i$.
\end{proof}

\begin{proposition}
We have $\mu=[\gamma,\U]$, where $\gamma \in \Gamma(M_0,\End(T^*M))$ is the
endomorphism of the tangent bundle of $M_0$ defined by
$$
\gamma dv^i = \half d\log\eta_i .
$$
\end{proposition}
\begin{proof}
By the definition of $\gamma$,
$$
\eta(\pi_i,\gamma^*\pi_j) = \half \eta_{ij} .
$$
We must show that
\begin{equation} \label{mu}
\eta(\pi_i,\mu^*\pi_j) = \eta(\U^*\pi_i,\gamma^*\pi_j) -
\eta(\gamma^*\pi_i,\U^*\pi_j) = \half (v^i - v^j) \eta_{ij} .
\end{equation}
We have
\begin{align*}
\eta(\pi_i,\mu^*\pi_j) &= (1-\tfrac{r}{2}) \eta(\pi_i,\pi_j) -
\eta(\pi_i,\nabla_jE) \\
&= (1-\tfrac{r}{2}) \delta_{ij} \eta_i -
\eta(\pi_i,[\pi_j,E]) - \eta(\pi_i,\nabla_E\pi_j) \\
&= - \tfrac{r}{2} \delta_{ij} \eta_i - \sum_k v^k \Gamma_{kji} \\
&= - \half \biggl( r \delta_{ij} \eta_i + \sum_k v^k
\bigl( \delta_{ik} \eta_{kj} + \delta_{ij} \eta_{jk} - \delta_{jk}
\eta_{ik} \bigr) \biggr) \\
&= - \half \bigl( \delta_{ij} \bigl( r\eta_i+E(\eta_i) \bigr) + (v^i-v^j)
\eta_{ij} \bigr) ,
\end{align*}
and \eqref{mu} follows.
\end{proof}

\subsection{The vector fields $\L_k$ on a semisimple Dubrovin manifold}

We now give a new proof of an important recent theorem of Dubrovin and
Zhang (\cite{DZ2}, Theorem 3.10.20). This proof is simpler than
theirs, though it is in the same spirit.
\begin{theorem} \label{unique}
  Let $M$ be a conformal semisimple Dubrovin manifold satisfying
  Condition \ref{gamma}. If $f\in\CO_\infty$ satisfies $\L_kf=0$ for
  all $k\ge-1$, then $f$ is a constant.
\end{theorem}
\begin{proof}
The theorem is proved by induction: if $f\in\CO_n$ satisfies $\L_kf=0$ for
all $k\ge-1$, then $f\in\CO_{n-1}$. In order to carry out this induction,
we must calculate the action of $\L_k$ on $\CO_n$ modulo $\CO_{n-1}$. By
\eqref{Lk}, we see that
\begin{align*}
\L_ku^a_n &= - \sum_{j=1}^n \p^{j-1} \* \Res_{z=0}
\bigl\<e,\delta_z^{k+1}(\p+z\X)^{n-j}\X du^a\bigr\> - \p^n \<
e,\U^{k+1}du^a \> \\
&= - n \Res_{z=0} \bigl\<e,\delta_z^{k+1}\*\p^{n-1}\X du^a\bigr\> - \p^n \<
e,\U^{k+1}du^a \> \pmod{\CO_{n-1}} \\
&= - n \sum_{\ell=0}^k
\bigl\<e,\U^\ell\bigl(\mu+\half\bigr)\U^{k-\ell}\p^{n-1}\X
du^a\bigr\> - \p^n \< e,\U^{k+1}du^a \> \pmod{\CO_{n-1}} .
\end{align*}
Form the generating function
\begin{align*}
\sum_{k=-1}^\infty \lambda^{-k-2} \L_k u^a_n
&= - n \, \bigl\< e , (\lambda-\U)^{-1}\bigl(\mu+\half\bigr) (\lambda-\U)^{-1}
\p^{n-1}\X du^a \bigr\> \\
& - \p^n \bigl\< e , (\lambda-\U)^{-1}du^a \bigr\> \pmod{\CO_{n-1}} .
\end{align*}

We now invoke the hypothesis that $M$ is semisimple. Since $\U
dv^i=v^i\,dv^i$, we see that
$$
\bigl\< e , (\lambda-\U)^{-1}du^a \bigr\> = \sum_{i,j} \bigl\< \pi_j ,
(\lambda-v^i)^{-1} J^a_i dv^i \bigr\> = \sum_i \frac{J^a_i}{\lambda-v^i} . 
$$
It follows that
$$
\p^n \bigl\< e , (\lambda-\U)^{-1}du^a \bigr\> = \sum_i 
\frac{v^i_n J_i^a}{(\lambda-v^i)^2} + \sum_i \frac{\p^nJ_i^a}{\lambda-v^i}
\pmod{\CO_{n-1}} .
$$
We also have
\begin{align*}
\p^{n-1}\X &= \p^{n-1} \bigl( u^c_1\AA_{bc}^a \bigr) \p_a\o du^b \\
&= u^c_n \AA_{bc}^a \, \p_a \o du^b = \sum_i v^i_n \, \pi_i \o dv^i
\pmod{\CO_{n-1}} ,
\end{align*}
hence
$$
\p^{n-1}\X du^a = \sum_i v^i_n J^a_i dv^i \pmod{\CO_{n-1}} .
$$
It follows that
\begin{multline*}
\bigl\< e, (\lambda-\U)^{-1}\bigl(\mu+\half\bigr)
(\lambda-\U)^{-1}\p^{n-1}\YY du^a \bigr\>
= \sum_{i,j} \frac{v^i_n J^a_i \bigl( \mu_j^i + \half \delta_j^i\bigr)}
{(\lambda-v^i)(\lambda-v^j)} \\
= \frac12 \sum_i \frac{v^i_n J^a_i}{(\lambda-v^i)^2} 
+ \sum_{i,j} v^i_n J^a_i\gamma_j^i \biggl( \frac{1}{\lambda-v^i} -
\frac{1}{\lambda-v^j} \biggr) \pmod{\CO_{n-1}} .
\end{multline*}
Assembling these results, we see that
$$
\sum_{k=-1}^\infty \lambda^{-k-2} \L_k u^a_n = - \bigl( \tfrac{n}{2} + 1
\bigr) \sum_i \frac{v^i_n J^a_i}{(\lambda-v^i)^2} - \sum_i
\frac{c^a_i}{\lambda-v^i} \pmod{\CO_{n-1}}
$$
for some coefficients $c_i^a\in\CO_n$.

Now suppose that $f\in\CO_n$ satisfies $\L_kf\in\CO_{n-1}$ for all
$k\ge0$. Extracting the coefficient of $(\lambda-v^i)^{-2}$ in
$$
\sum_{k=-1}^\infty \lambda^{-k-2} \L_kf = \p_{n,a}f \sum_{k=-1}^\infty
\lambda^{-k-2} \L_ku^a_n \pmod{\CO_{n-1}} ,
$$
we see that $\bigl( \frac{n}{2} + 1 \bigr) v^i_n J^a_i\,\p_{n,a}f=0$, and
hence that $\p_nf=0$. In other words, $f\in\CO_{n-1}$.
\end{proof}

\begin{corollary}
On a semisimple Frobenius manifold, any solution $\{\PHI_g\in\CO_\infty\mid
g\ge1\}$ of the Virasoro constraints (Theorem \ref{Virasoro}) and the
dilaton equation \eqref{Dilaton} is unique (except in genus $1$, where it
is determined up to an additive constant).
\end{corollary}
\begin{proof}
  Consider two solutions $\{\PHI_g\}$ and $\{\tilde{\PHI}_g\}$ of the
  Virasoro constraints; we argue by induction on $g$ that
  $\PHI_g=\tilde{\PHI_g}$. Suppose that $\PHI_h=\tilde{\PHI_h}$ for
  $h<g$. We see that $\L_k(\PHI_g-\tilde{\PHI_g})=0$ for all $k\ge-1$,
  hence by Theorem \ref{unique}, $\PHI_g-\tilde{\PHI}_g$ is a
  constant.  If $g>1$, the constant vanishes by the dilaton equation,
  which implies that
\begin{equation*}
\Dil(\PHI_g-\tilde{\PHI}_g) = (2g-2)(\PHI_g-\tilde{\PHI}_g) .
\qedhere\end{equation*}
\end{proof}

\subsection{The function $\GG$ on a semisimple Frobenius manifold}

A solution of the genus $1$ Virasoro constraints \eqref{V1} on a
semisimple Frobenius manifold was found by Dubrovin and Zhang
\cite{DZ}; by Theorem \ref{unique}, their solution is the unique one.
Their formula for $\GG$ involves a function on the semisimple locus
$M_0$ called the isomonodromic $\tau$-function $\tau_I$; in defining
it, we follow the exposition of Hertling \cite{Hertling}.
\begin{lemma}
The differential form
$$
\alpha = \tfrac{1}{8} \sum_{i,j} (v^i-v^j) \frac{\eta_{ij}^2}{\eta_i\eta_j}
dv^i \in \Om^1(M_0)
$$
is closed, and $\alpha(X) = - \quart \Tr(\mu[\gamma,\A{X}])$
\end{lemma}
\begin{proof}
We have
\begin{align*}
d\alpha &= \tfrac{1}{4} \sum_{i,j,k} \frac{\eta_{ij}}{\eta_i\eta_j} \biggl(
\eta_{ijk} - \half \frac{\eta_{ij}\eta_{ik}}{\eta_i} - \half
\frac{\eta_{ij}\eta_{jk}}{\eta_j} \biggr) \, (v^i-v^j) \, dv^k \. dv^i \\
&= \tfrac{1}{8} \sum_{i,j,k}
\frac{\eta_{ij}\eta_{jk}\eta_{ki}}{\eta_i\eta_j\eta_k} \, (v^i-v^j) \, dv^k
\. dv^i .
\end{align*}
This vanishes, since $(v^i-v^j)\,dv^k\.dv^i$ is antisymmetric in the
indices $i$, $j$ and $k$.

The matrix $\A{i}$ has components
$(\A{i})_j^k=\delta_{ij}\delta_i^k$. Since the matrices $\U$ and $\A{i}$
are diagonal in the frame $dv^i$ of $T^*M_0$, we have
\begin{align*}
\p_i\log\tau_I &= - \quart \Tr([\gamma,\U][\gamma,\A{i}]) \\ &= -
\quart \Tr( \gamma\U\gamma \A{i} - \gamma\U \A{i}\gamma - \U\gamma^2
\A{i} + \U\gamma \A{i}\gamma ) \\
&= - \quart \sum_j \bigl( \gamma_j^i\U_j^j\gamma_i^j -
\gamma_j^i\U_i^i\gamma_i^j - \U_i^i\gamma_i^j\gamma_j^i +
\U_j^j\gamma_j^i\gamma_i^j \bigr) \\
&= \half \sum_j ( v^i - v^j ) \gamma_j^i\gamma_i^j = \<\pi_i,\alpha\> .
\qedhere\end{align*}
\end{proof}

The \textbf{isomonodromic tau-function} $\tau_I$ is defined by the
equation
$$
d\log\tau_I=\alpha .
$$
Note that just like $\GG$, $\log\tau_I$ is only determined up to an
additive constant.
\begin{theorem}
$\GG = \log\tau_I - \tfrac{1}{48} \sum_i \log\eta_i$
\end{theorem}
\begin{proof}
We calculate the derivative of the two terms contributing to $\L_k\GG$
separately: we will see that they respectively contribute the two terms of
\eqref{V1}. We have
\begin{align*}
\L_k\log\tau_I &= \quart \Tr(\mu[\gamma,\A{E^{\circ(k+1)}}])
= \quart \Tr(\mu[\gamma,\U^{k+1}]) \\
&= \quart \sum_{\ell=0}^k \Tr(\mu\U^\ell[\gamma,\U]\U^{k-\ell}) 
= \quart \sum_{\ell=0}^k \Tr(\mu\U^\ell\mu\U^{k-\ell}) .
\end{align*}
Further, since $e(\eta_i)=0$, we have
\begin{align*}
\half \L_k(\log\eta_i) &= \half \bigl( \L_k+(v^i)^{k+1} e \bigr)
(\log\eta_i)
= \sum_j \bigl( (v^i)^{k+1} - (v^j)^{k+1} ) \gamma_j^i \\
&= \sum_j \frac{(v^i)^{k+1} - (v^j)^{k+1}}{v^i-v^j} \mu_j^i
= \sum_{\ell=0}^k \sum_j (v^j)^\ell \mu_j^i (v^i)^{k-\ell} \\
&= \sum_{\ell=0}^k \< e , \U^\ell \mu \U^{k-\ell}dv^i \> .
\end{align*}
Summing over $i$, the Virasoro constraints $z_{k,1}$ follows, since
$\om=\sum_i dv^i$.
\end{proof}

It is proved in \cite{DZ0} that $\GG$ automatically satisfies the
differential equation $\Psi\GG=0$ of \eqref{elliptic}. Ideally, this
should be part of a larger phenomenon, whereby all of the known
differential equations satisfied by the Gromov-Witten potentials
$\{\PHI_g^X\}_{g\ge1}$, namely the Virasoro constraints and the
topological recursion relations, have a (necessarily unique) solution
on every semisimple Frobenius manifold.

\end{document}